\documentclass[12pt,onecolumn]{IEEEtran}
\usepackage{atbegshi}
\AtBeginDocument{\AtBeginShipoutNext{\AtBeginShipoutDiscard}}
\addtocounter{page}{-1}
\usepackage[american]{babel}
\usepackage[utf8x]{inputenc}
\usepackage{amsmath}
\usepackage{graphicx}
\usepackage[colorinlistoftodos]{todonotes}
\usepackage{amsmath,amssymb,euscript,yfonts,psfrag,latexsym,dsfont,graphicx,bbm,color,amstext,wasysym,pdfsync,soul,framed,url,balance,soul}
\usepackage{apacite}
\usepackage{lipsum}
\usepackage{csquotes}
\usepackage{soul}

\def\R{\mathbb{R}}

\newcommand{\beq}{\begin{equation}}
\newcommand{\eeq}{\end{equation}}
\newcommand{\bea}{\begin{eqnarray}}
\newcommand{\eea}{\end{eqnarray}}

\newcommand{\bsea}{\begin{subeqnarray}}
\newcommand{\esea}{\end{subeqnarray}}

\newcommand{\anne}{\color{blue}}

\def\bmat{\left[ \begin{array}}
\def\emat{\end{array} \right]}

\newcounter{acount}

\begin{document}
\title{Harmonic Analysis of Social Cognition}
\author{Anne Maass, Michele Pavon and Caterina Suitner}
\thanks{A. Maass and C. Suitner are with the 
Department of Developmental Psychology and Socialisation, University of Padova, Padova, Italy. {\tt\small  anne.maass@unipd.it, caterina.suitner@unipd.it}.}\thanks{M. Pavon is with
the Department of Mathematics ``Tullio Levi-Civita'', University of Padova, Padova, Italy. {\tt\small
pavon@math.unipd.it}}
\maketitle
\begin{abstract}In this paper, we argue that  some fundamental concepts and tools of {\em signal processing} may be effectively applied to represent and interpret 
{\em social cognition processes}. From this viewpoint, individuals or, more generally, social stimuli are thought of as a weighted sum of harmonics with different {\em frequencies}: Low frequencies represent general categories such as gender, ethnic group, nationality, etc., whereas high frequencies account for personal characteristics. Individuals are then seen by observers  as the output of a {\em filter} that emphasizes a certain range of high or low frequencies. The selection of the filter depends on the {\em social distance} 
between the observing individual or group and the person being observed  as well as on {\em motivation}, {\em cognitive resources} and {\em cultural background}. 
Enhancing low- or high-frequency harmonics is not on equal footing, the latter requiring  {\em supplementary energy}. This mirrors a well-known property of signal processing filters. More generally, in the light of this correspondence, we show that several established results of social cognition admit a natural interpretation and integration in the signal processing language. While the potential of this connection between an area of social psychology and one of information engineering appears considerable
(compression, information retrieval, filtering, feedback
, feedforward, sampling, aliasing, etc.), in this paper we shall  limit ourselves to laying down what we consider the pillars of this bridge on which future research may be founded. 
\end{abstract}

\begin{IEEEkeywords} Signal processing, social cognition, filtering, person perception, construal level theory.
\end{IEEEkeywords}
\section{Introduction}\label{introduction}
Social environments are extremely complex, yet people manage to successfully navigate them by selecting or filtering out information according to the relevance to their current aims and by flexibly shifting attention from specific details to the larger picture, or vice versa. How people prevent information overload and how they select and process social signals at different levels of abstraction remains an important challenge for social cognition. In this article, we propose signal processing theory as an overarching framework, and as a common language, to integrate different lines of social cognition that have so far been studied in relative isolation. The main idea is that  information that we actively seek or incidentally encounter is carried by signals that we select, process, archive, and transmit. If well managed, complex signals are harmonically orchestrated in a comprehensible and efficient way that provides information; otherwise they are noise. We here propose an application of the main principles of harmonic analysis to the management of social signals, arguing for the need to integrate knowledge from disciplines that are traditionally conceived as distant, but that may greatly gain from reciprocal inputs, both under a theoretical and practical perspective. For simplicity we will refer to this approach as social signal-processing, henceforth SSP.

We are by no means the first to apply signal-processing to psychological phenomena. One of the most notable precursors is  signal detection theory \cite{swets2014signal}, which has been applied to different psychological realms, including memory, attention, decision making, eyewitness identification.  Signal processing approaches  have also been used in human vision \cite{schyns1998development,schyns1999dr,pylyshyn1999vision,serre2007feedforward,meese2009tutorial,doutsi2018retina}. A third important line of applications lays at the interface between psychology and engineering, such as artificial social intelligence \cite{vinciarelli2009social}\footnote{Notice that the latter research deals with the {\em automatic} processing of social signals such as turn taking, 
smiles, blinks, posture, etc., whereas we apply signal-processing theory to {\em human} social cognition.}. However, to the best of our knowledge, the theory and language of signal processing has never been applied to social cognition.

The paper is organized as follows. In three distinct sections, we introduce a non-technical overview of signal processing theory together with the proposed applications  to social signals. We introduce, step by step,  the main tenets of the theory, each followed by a selective overview of social-cognitive phenomena that, in our opinion, are interpretable within the signal processing framework. To illustrate the utility of the integrative framework, we present, for each tenet,  different subareas of social cognition that can be interpreted as reflecting the same signal processing principles. We start from the idea that human beings can be envisaged as a sum of harmonics of different weight and frequency (see Section \ref{representation}, entitled Representation). We  introduce the reader to some basic concepts of signal processing, to then argue that social targets can be viewed as a multitude of harmonics whose relative weight varies as a function of intrinsic variables and of the social context in which the target is embedded. In Section \ref{filtering}, we then shift attention from the target to the observer. We introduce the concept of {\em filtering}  and  argue that observers systematically filter social signals in a way that harmonics in a certain frequency band  are enhanced, whereas others are  attenuated or even completely  inhibited. In Section \ref{what filter and why}, entitled “What Filter and Why”, we introduce the idea that the selection of the filter employed by the observer is affected by (social) distance, motivation and cognitive resources, and culture  with energy playing a crucial role in the selection. Indeed, we advance the hypothesis that low-pass filters, concentrating on general categories, are {\em passive}, namely they do not require further energy supply. Instead, high-pass filters, concentrating on specific characteristics, are  {\em active}, i.e. require a supplementary quantity of energy.  In the final Section \ref{conclusions}, we first discuss the similarities and differences between SSP and Construal Level Theory \cite{trope2010construal}. We point out that, although there are important points of contact, neither theory includes the other, and, in some cases, the two models produce distinct predictions. We then derive seven novel hypotheses from SSP, which can be investigated in future research.




Throughout this paper, we  keep the technical details to a bare minimum, while referring the interested reader to  the Appendix for a more detailed descriptions of the  concepts and models of signal processing. Also, the relevant social-cognitive literature is reviewed in a highly selective fashion, without claiming completeness or representativeness. Besides reporting on recent work, we also cite classical research from the 80s and 90s, a period that may be considered the “golden age” of social cognition, when impression formation and categorization vs. individuation were at the center of the field. The reported examples are meant to illustrate signal processing interpretations of diverse phenomena while keeping in mind both the potential and the limits of our unifying framework. 
Throughout this article, we pursue two main aims: On the one side, we want to offer a comprehensive framework and language that can accommodate diverse socio-cognitive phenomena lacking theoretical integration. On the other side, we want to advance novel, testable hypotheses and methodological and analytical techniques that are suggested by the analogy with signal processing theory.

\section{Representation}\label{representation}
\subsection{Signals}\label{signals}

A {\em signal} is a map from  a time or spatial domain taking values in a finite or infinite set. For instance, the human voice, an electrocardiogram, the level of tide in Venice, the daily count of new positive Covid-19 cases, the monthly unemployment rate, a binary sequence, a picture, a video are all signals. The word signal is preferred in information engineering to the more general {\em map} or {\em function} terms to emphasize the fact that these maps {\em carry information}.  In this paper, we advance the idea that social stimuli may also be viewed as signals.

Today's world is dominated by the necessity to {\em elaborate} signals. Signals are transformed by systems in order to amplify them, to transduce them, to remove or attenuate noise, to make them suitable for storage or transmission (e.g., compression), to encrypt them for security purposes, etc. Indeed, most of the devices that are so crucial in our daily  life such as cell phones, computers, appliances, cars, microphones, cameras, audio sets,  etc. feature simple to rather  complex on-line or off-line elaboration of electrical, mechanical, chemical, hydraulic, etc. signals. The same does, in a marvelous way, our body elaborating signals so as to keep a nearly constant temperature in all conditions, to control metabolic processes, to retrieve information and to send new one  to our brain, to name but a few. Coding and processing of signals have also characterized the evolution of human communication, starting from tuning gestural and vocal signals to exchange information with our fellows \cite{tomasello2010origins}. We will argue below in Subsection III \ref{second hypothesis} that social inputs are also  "elaborated" in that the observer processes the social stimulus so as to emphasize certain aspects of it and attenuate or even remove others.


\subsection{Harmonics}\label{harmonics}
The first step in Signal Processing  to elaborate a signal consists in {\em representing} the signal as a finite or infinite weighted sum, called {\em Fourier series,} of simpler ({\em canonical}) signals\footnote{This idea originates in the 18th century with Euler, D'Alembert and D. Bernoulli. Fourier Series are discussed in the Appendix A.2, see, in particular, formula (\ref{FSeries}).}. These canonical signals are called {\em harmonics}, a term borrowed from music. The fundamental reason for doing so is that canonical signals are transformed by a large class of systems  in a much simpler way than arbitrary signals. A system $\Sigma$ 
\begin{center}
\begin{equation}\label{System}
\setlength{\unitlength}{.2cm}
\begin{picture}(43,8)(0,5)
\put(16.5,10){\makebox(0,0){$x(t)$}}
\put(31,10){\makebox(0,0){$y(t)$}}
\put(12.5,8){\vector(1,0){7.5}}
\put(28,8){\vector(1,0){7.5}}
\put(20,6){\framebox(8,4){$\Sigma$}}
\end{picture}
\end{equation}
\end{center} 
is any transformation of an input signal $x(t)$ into an output signal $y(t)$. From the transform of the canonical signals, if the system satisfies the {\em superposition principle}, i.e. it is linear (see the Appendix \ref{convolutional}), it is then possible to compute the transform of the original signal.   

\begin{figure}
\begin{center}
	\includegraphics[scale=1.5]{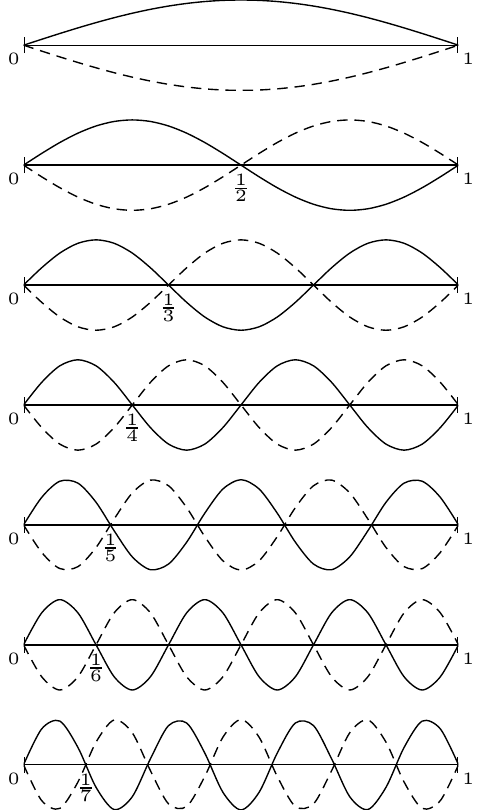}
	\caption{Normal modes of a vibrating string.}
	\label{vibratingstring}
	\end{center}
\end{figure}

\begin{figure}
\begin{center}
	\includegraphics[scale=2]{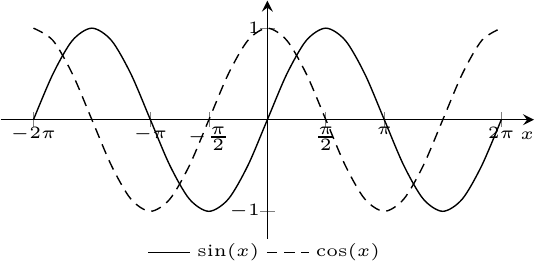}
	\caption{Sine and cosine signals}
	\label{sine_cosine}
	\end{center}
\end{figure}

The first historical example of harmonics was provided by Daniel Bernoulli in 1753. He stated: ``The general motion of a vibrating system is given by a superposition of its proper vibrations”.  For a vibrating string, clamped at both ends, the first {\em normal modes} (proper vibrations) for $k=\pm1,\pm2,\pm3,\pm4,\pm5,\pm6,\pm7$ are depicted in Figure \ref{vibratingstring} (the normal mode corresponding to $k=0$ is the horizontal string).
Mathematically, normal modes are the sinusoidal signals\footnote{The oscillation of a frictionless  spring-mass system around  equilibrium is also  described by a sine wave.}

\[x_k(t)=\sin\left( k\pi t\right), \quad k=0,\pm1,\pm2,\pm3,\ldots
\]
where $\sin(t)$ and $\cos (t)=\sin(t+\frac{\pi}{2})$ are depicted in Figure \ref{sine_cosine}. 
Notice that, in the case of the vibrating string, $t$ varying in $[0,1]$ (Figure \ref{vibratingstring}) is not time but a {\em spatial} variable.
If $k>0$, the frequency of $x_k$ and of $x_{-k}$ is  precisely $k$. For instance, $x_4$ and $x_{-4}$ have frequency which is four times the frequency (one) of $x_1$ and $x_{-1}$ which are called {\em first harmonics}. {\em Low-frequency harmonics} have small absolute value of $k$, denoted $|k|$, (see Appendix \ref{real}). High-frequency harmonics have  large $|k|$. 
This means that any “wave”   $\{x(t), 0\le t\le 1\}$ described by the string can be expressed as a  weighted sum of the fundamental waves corresponding to the normal modes (Fourier series) of the form
\begin{equation}\label{stringFS}x(t)=\sum_k a_k x_k(t), \quad 0\le t \le 1.
\end{equation}
In this simple classical example, we already see one crucial difference between low-frequency and high-frequency harmonics. On the same one-dimensional spatial interval $[0,1]$, low-frequency harmonic vary {\em slowly} whereas high-frequency harmonics vary {\em rapidly}. For instance, the first harmonic has only one maximum or minimum. The seventh has seven between maxima and minima. This difference can also be described as follows: The low-frequency harmonics are {\em smoother}  than  the high-frequency ones. This is captured by the fact that the slope of the tangent to the graph (derivative of the harmonics) takes larger (absolute) values in the case of a high-frequency harmonics. This is signaled by the curve increasing or decreasing rapidly.

To  start understanding the analogy with music, consider a mixed choir consisting of several parts corresponding to the different vocal ranges.  Sometimes the entire choir sings the main melody producing monophonic music where all notes sung are in {\em unison}. The latter means that the different parts of the choir  are either singing the same pitch or {\em pitches that are  one or more octaves} apart. An octave is the interval  between one musical pitch and another with {\em double its frequency}\footnote{The frequency of oscillations in a wave motion is measured in Hertz. The frequency is $1$ Hz when one oscillation occurs in one second. When $1000$ oscillations occur in one second, the frequency is $1000$ Hz, or $1$ kHz (kilohertz).}. The human ear tends to hear  two notes  separated by one octave as if they were equal. This is why we label in the same way  notes that are exactly one octave apart. Thus, when listening to such a choir, we hear a tone, which is the superposition  of pure tones\footnote{A sound that includes only one frequency.} that are an integer number of octaves apart. The weight on each harmonic can be thought of as the sound volume of that particular harmonic  measured in decibels. 

Consider now  the following example: There are currently thousands of  languages in the world that may be considered separately, representing the high-frequency range. However, based on their similarities and shared origins, these same languages may also be grouped into Italic, Celtic, Slavic etc. languages (intermediate frequency range), which, in turn may be categorized in broad language families, such as Indo-European, Niger-Congo, Austronesian, Trans-New Guinea, Sino-Tibetan, and Afro-Asiatic languages (low-frequency).

Keeping in mind Figure \ref{vibratingstring}, let us return for a moment  to music. In  a stringed instrument like a guitar, each cord has an associated fundamental frequency and corresponding first harmonic. We hear the latter if we pluck the string.  Higher frequency harmonics are played by touching (but not fully pressing down the string) at an exact point on the string while plucking it: the pitch so produced is always higher than the fundamental frequency of the string.

A system which accomplishes the complete or partial suppression of some harmonic components 
of the input signal $\{x(t)\}$ is called {\em filter}. 
A very important class of filters (convolutional systems, see Appendix \ref{convolutional}) can change {\em selectively} the energy of the single harmonic by replacing each coefficient $a_k$ with a coefficient $b_k$, see (\ref{selective}) in Appendix \ref{convolutional}. For instance, in an audio system, turning a suitable knob, we can increase the relative energy of the low-frequency harmonics (bass) or of the high-frequency harmonics (treble). This is accomplished by passing the audio signal through a suitable filter.

 We can use Fourier series to represent basically any signal defined on a bounded domain. While in some applications it might be apparent what the different frequency harmonics represent (vibrating string, vocal signal, etc.), this may not be so salient in others. For instance, an image can be represented through a two-dimensional Fourier series (\ref{2DFS}). High-frequency harmonics in the two directions account for sharp transitions while low-frequency ones represent smooth changes in the image, see  Subsection IV\ref{analogy}. This is all we need when applying signal processing to human  or  computer vision. Indeed, this abstract representation has been very successful  to describe the human visual system as well as in computational vision \cite{sonka2014image}. Throughout this paper, we shall take a similar standpoint when representing  target persons and social stimuli as a Fourier series of harmonics\footnote{See in the Appendix the end of A.2 for a discussion of a more realistic harmonic analysis model for social stimuli.}.

\noindent
{\em Summary: The first step to process a signal consists in representing the signal as a weighted sum (superposition) of simpler signals (harmonics) that oscillate at different frequencies. }

\subsection{The first hypothesis: Human beings  as a weighted sum of harmonics}\label{first hypothesis}
{\em Human beings, in the role of targets of social cognition, may be considered as a weighted sum of harmonics. Low-frequency harmonics represent general categories such as gender, nationality, ethnic group or religion. High-frequency harmonics represent specific behaviors, psychological states, or elements of appearance (e.g., clothing, hair style). Relatively stable individual characteristics or traits --- such as height, occupation, extraversion, temper, attachment style, or intelligence --- occupy the middle range. This representation is postulated, as the full range of all harmonics of a specific person may  never be accessible, partially because some information is unknown, partially because it would be so complex to result in noise. In the case of the vibrating string or an image, harmonics vary as a function of {\em spatial} variables. Here, in the postulated representation of social stimuli, harmonics vary, as in the case of a voice or music recording, as a function of {\em time}}.

In signal processing terms, human beings can be envisaged as the superposition of a myriad of harmonics of different frequency. At the high-frequency end are ephemeral elements such as transitory emotions  or behaviors performed at any specific moment and context. As frequency decreases, we are likely to find individualizing characteristics that are more stable over time than behaviors but that vary greatly from individual to individual. For instance, people have unique personalities, habits, traits, and possessions that reflect their individual predispositions and life histories and that distinguish them from other human beings.  On the low-frequency extreme of the continuum\footnote{Mathematically, of the sequence.}, we find characteristics that are common to many individuals and that define our social memberships. Examples of such low-frequency harmonics are categories such as nationality, faith, gender, or ethnic group, which are constant or vary very slowly over time. This continuum  from high to low-frequency harmonics is conceptually similar to hierarchical models in social cognition that conceive behaviors, traits, and categories as increasingly general and abstract types of information. In fact, hierarchies from the concrete to the abstract, or from the individualizing to the categorical, are built into numerous social-psychological theories, among which impression formation models  \cite{rothbart1981memory,fiske1990continuum}, social identity and self-categorization theories \cite{tajfel1979integrative, turner2011self}, tripartite models of selfhood \cite{sedikides2011individual}, construal level theory \cite{trope2010construal}, and the linguistic category model {\cite{semin1991linguistic}. According to SSP, across all of these models, the concrete or individualizing end of the continuum can be envisaged as high-frequency and the abstract or categorical end of the continuum as low-frequency harmonics.

At any given moment, all harmonics are theoretically available to the observer, but their relative weight will, in part, determine which harmonics will meet the observer's eye and which harmonics will go unnoticed. We argue that the relative weight of high versus low-frequency harmonics depends (a) on the properties of the signals to be processed, (b) on the target person's identity management, and (c) on the context in which the signal is embedded. We  focus here  on the properties of the stimulus and of the context in which it is encountered. The role of the observer in selecting and interpreting signals is discussed in  Subsection III\ref{second hypothesis}.

As far as the properties of social signals are concerned, some signals carry intrinsically more weight than others. For instance, extreme, negative and non-normative information captures attention more and has a stronger impact than moderate, positive or normative information (for an overview see \cite{alves2017good}). Many authors locate the reasons of this robust bias in the observer for whom negative information is particularly relevant,  an issue that we will come back to in  Subsubsection IV\ref{motivation}. In contrast, Alves and collaborators attribute it to the intrinsic properties of the stimulus, such that negative instances tend to be both, more extreme and more distinct, one being dissimilar from the other, whereas positive instances look alike and typically lie in the middle range, denoting moderation\footnote{ \cite{alves2017good} have this beautiful quote from Tolstoy's Anna Karenina: “Happy families are all alike;
every unhappy family is unhappy in its own way”.}. In a nutshell, “good is more alike than bad”  \cite[p.1171]{koch2016general} which echoes the fact that the low-frequency normal modes of the vibrating string vary slowly whereas the high-frequency ones vary rapidly, see  Subsection II\ref{harmonics}. For instance, attractive (vs. unattractive) faces, likeable (vs. unlikeable) groups, and positive (vs. negative) words tend to be more similar to each other  and tend to lie in the middle range of relevant dimensions, which correspond to the low-frequency range. In line with this argument, positive categories also tend to be broader than negative categories, as shown by studies in which participants sorted representative samples of positive and negative stimuli into groups \cite{koch2016general}. In SSP terms, this suggests that positive stimuli are more homogeneous and hence of lower frequency than negative stimuli whose unique and heterogeneous properties locate them in the high-frequency range. By extension, the relative weight of each harmonic is, to some degree, determined by the properties of the social signals (in our example: their valence) that are being processed. 

The weight of specific harmonics may also be influenced by the target person themselves. Human beings are actors on the social stage who actively engage in {\em identity management}. They may consciously or subconsciously stress their individual characteristics, by displaying their possessions or personal preferences, or they may signal their category memberships. They may present themselves as unique individuals or as more or less interchangeable parts of a collective. For instance, bumper stickers such as “low standards good life”  reveal relatively high-frequency information, whereas “sorry, girl, I am gay” or “Black owned business” underline low-frequency category memberships. People may choose to disclose or conceal  “invisible”  categories such as religion, sexual orientation, mental illness \cite{newheiser2014hidden}.  Importantly, self-construals are quite fluid and variable, depending on context \cite{cadinu2013chameleonic}, personality e.g., the tendency to monitor one's impression \cite{snyder1987public}), and culture \cite{markus1991culture}. This and related research suggests that people may strategically communicate different aspects of their identity, giving greater weight either to the qualities that make them unique or to those shared with others. 

Finally, relative weight is determined by the {\em context} in which the social signals are embedded. Any kind of social signal unfolds in a specific information environment that affects the weight of the signal. Three examples are {\em solo} status, category accentuation, and group context. There is ample evidence that, as proportionate group size becomes unbalanced, low-frequency harmonics (social category membership) become prominent, obscuring individuating information. For instance, gender or race become highly salient in numerically unbalanced groups or under {\em solo} status (e.g., one woman interacting with 5 men, or one White individual interacting with 5 Black individuals). {\em Solo} status increases the likelihood that observers will rely on categorical rather than individuating information and, as a consequence, will stereotype the target as a member of that category, see \cite{taylor1978categorical,taylor1981categorization}. {\em Solo} status will also shift the target’s self-perception, such that targets will align their self-construal with the category they are members of, e.g., \cite{niemann1998relationship}, which may potentially have tangible consequences for their performance \cite{sekaquaptewa2003solo}. In both cases, the relative weight of low- versus high-frequency harmonics is a direct function of the social (here: numerical) context. 

A different way in which context may shift the weight from high to low-frequency is through verbal labels or visual boundaries that create larger categories. A classic example is Tajfel and Wilkes’ groundbreaking study, in which 8 lines of steadily increasing length where labeled as A or B either in a random fashion or systematically, so that the shorter lines were labeled A and the longer lines B \cite{tajfel1963classification}. In the latter condition, the authors created two artificial categories, which led people to exaggerate the actual difference between the short vs. long lines, a phenomenon known as category accentuation. Subsequent studies in which arbitrary boundaries were cued either through verbal labels or though visual boundaries \cite{rothbart1997effects,foroni2011category} found not only inter-category accentuation, but also intra-category minimization. In our SSP framework, such arbitrary boundaries shift the weight from the high to the low-frequency end of the range.

A third example derives from a recent set of studies by Cocley and Payne \cite{cooley2019group}, showing that individuals are judged more in line with category stereotypes when part of a group rather than when shown alone. For instance, the same Black individuals were perceived as more typical of their category when shown in the context of other Black individuals than when shown alone. This suggests that the simple presence of other members of the same category shift the weight from high towards low-frequency harmonics.

Together, the above findings suggest that there are different ways in which high or low-frequencies of the social stimulus become salient. These include the intrinsic properties of the stimulus (such as the case of negative vs. positive stimuli), the target’s impression management emphasizing individualizing or categorical features, or the context (such as solo status, arbitrary boundaries) that changes the relative salience of individualizing or category information. These variations in stimulus properties, impression formation, and context can jointly be interpreted as tuning the weights with which harmonics appear in the representation of an individual or, more generally, of social stimuli. 

So far, we have focused on the features (harmonics) of the social target, whereas in the next section we will analyze how observers select certain types of information when elaborating a social stimulus, a mechanism called filtering in signal processing. 

\section{Filtering}\label{filtering}

\subsection{Frequency selective and frequency shaping filters}\label{frequency selective}

Systems (\ref{System}) are classified according to their properties. We only mention here: Stability, Causality, and Passivity.  BIBO-stability amounts to the property that bounded inputs $\{x(t)\}$ are mapped into bounded outputs $\{y(t)\}$. 
 BIBO-stability is a fundamental requirement for any designed system where all signals, whatever their nature, must remain within bounds to avoid damage/failure. For convolutional systems (\ref{conv}), see Appendix \ref{convolutional}, BIBO-stability implies a form of input-output continuity: Small perturbations of the input cause small perturbations of the output (for the role of this  property  in the social cognition literature see Subsection IV\ref{third and fourth} below). Causality is the property that, at each time $t$, the output $y(t)$ depends only on the values of the input prior to time $t$, namely $\{x(s), s\le t\}$. Finally, a system is called {\em passive} if it does not require a supplementary energy source (see Subsection IV\ref{third and fourth} below for the role of this property in social cognition). Filters are special systems which may or may not possess the above properties.

Consider first the so-called {\em ideal} (frequency selective) filters which completely suppress harmonics with frequencies in certain bands ({\em stopband}), passing the others unaltered in energy ({\em passband}). The magnitude of the frequency response (cf. (\ref{FR}) in Appendix \ref{convolutional})  of some fundamental ideal filters is depicted in Figure \ref{IF} as a function of the positive circular frequency $\omega=2\pi f$. The graph of such magnitude is usually represented, like here, only for positive frequencies as it is symmetric with respect to the vertical axis.
\begin{figure}
\begin{center}
	\includegraphics[scale=2]{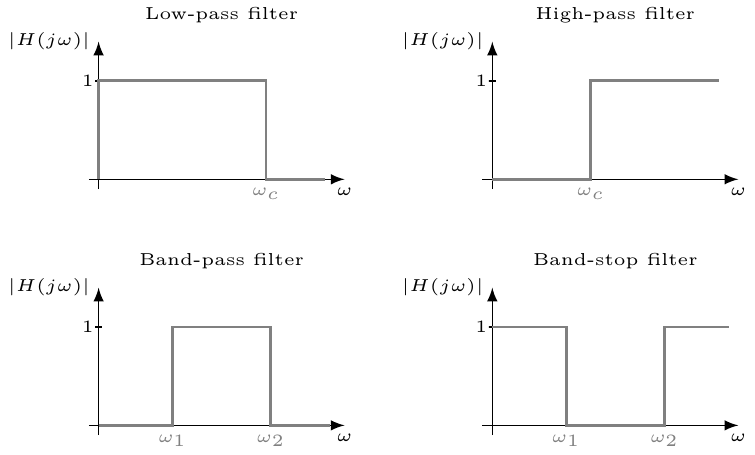}
	\caption{Ideal Filters: Frequency response magnitude}
	\label{IF}
	\end{center}
\end{figure} 
Ideal filters only serve as a utopian reference as far as their frequency response magnitude goes. Indeed, while they cannot be physically realized, they would also have undesirable features such as being unstable and non causal, cf. \cite[Section 3.9]{oppenheim1997signals}. Usually, non-ideal filters that are stable and causal and whose frequency response approximates that of an ideal filter are used instead.
{\em Damped harmonic oscillators}, for instance, see Appendix \ref{low vs. high filters}, are a class of systems which approximate  ideal low-pass filters. They abound in nature and in human-made devices such as the damped mass-spring system described in Appendix \ref{low vs. high filters}. They are stable, causal and passive.

Consider now, as high-pass filter, the ideal differentiator
\[y(t)=\frac{dx}{dt}(t)
\] 
which associates to each smooth signal its derivative (slope of the tangent to the graph). This filter  is not BIBO-stable but, in some mitigated form (especially in the digital version), it is used in control systems because of certain positive features such as making systems with a lot of inertia more prompt to respond. The ideal differentiator may be seen (Appendix \ref{low vs. high filters}) to change the energy of the  harmonics {\em proportionally to their frequency}: the higher the frequency of the harmonic the higher the increase of energy. It can, therefore, be used in image processing, to identify objects in images through {\em edge detection}: In a black and white picture, it helps localizing large local changes in the grey level image which correspond to the edges. While a rough approximation of a differentiator at the low frequencies may be realized through a simple, passive circuit, cf. (\ref{RC}), a better approximation requires
using an active element such as an {\em operational amplifier}, see Figure \ref{active}.
\begin{figure}
\begin{center}
	\includegraphics[scale=1.5]{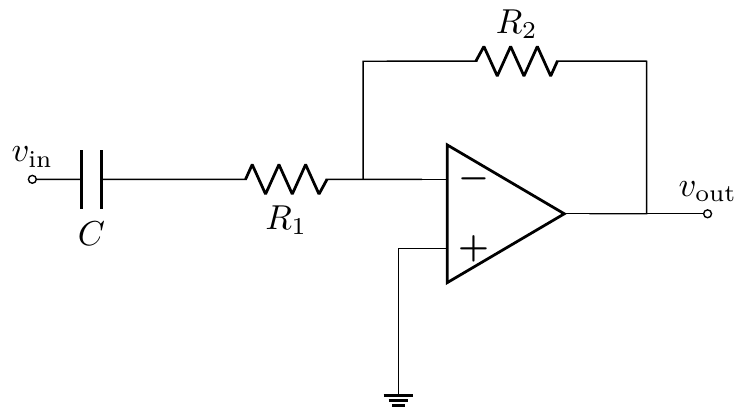}
	\caption{An active high-pass filter.}
	\label{active}
	\end{center}
\end{figure}
Band-pass and band-stop filters, whose magnitude responses are depicted in the second row of Figure \ref{IF}, also play a significant role in SSP. They will be introduced and discussed directly in connection with social cognition processes starting form  Subsubsection III \ref{beyond}. 

\noindent
{\em Summary: Filters are systems transforming signals so as to pass unaltered or enhance harmonics in a certain frequency band (passband) while attenuating or even deleting them in other frequency bands (stopband). Low-pass filters induced by a damped oscillator (or an RLC circuit, see  Appendix \ref{low vs. high filters}.) abound in nature, are easy to realize and are intrinsically stable (a small perturbation of the input causes a small perturbation of the output). They are, moreover, {\em passive}, i.e. not requiring an external power \footnote{Energy rate.} source. High-pass filters may be unstable or cause instability of an interconnection\footnote{Although there exist unstable low-pass and stable high-pass filters, stability is much easier to attain when constructing a low-pass rather than a high-pass.}. Moreover, a precise realization of an ideal differentiator involves introducing {\em active} elements. These require an external power source. Thus there is a clear asymmetry between the low- and high-pass filters described in this subsection.}

\subsection{The second  hypothesis: Observers and filters}\label{second hypothesis}

{\em Observers do not perceive targets as the full complexity of the weighted sum of their harmonics as $x(t)$ in (\ref{stringFS}), but rather as the output $y(t)$ of a filter fed by the harmonic sum as  in (\ref{System}). The filter enhances harmonics in a certain frequency band and  simultaneously attenuates them  or even cancels them in others.}

So far, we have focused on the stimulus person, described as  a weighted sum of harmonics, and on factors that (objectively) shift the relative salience of high- vs. low- frequency harmonics. We  shall now direct attention to the observer who spontaneously or intentionally selects specific frequency ranges while neglecting others. Observers are unlikely to ever have access to all harmonics of a target person, given that each person possesses an almost infinite number of features, many of which invisible. Yet, to form an impression, observers are likely to draw on a wide range of data and select information in a specific frequency range.
The idea that observers filter the incoming information is by no means new to social cognition. Social-cognitive research and theorizing has in large parts focused on the perceiver rather than on the target and has stressed the fact that observers are are “cognitive misers” \cite{fiske2013social} who simplify their tasks by using mental shortcuts, among which the selective focus on certain stimuli rather than on others. Selective filtering of specific stimulus ranges is at the heart of different socio-cognitive theories as we will illustrate below. We will cite four instances of such band-specific filtering, stemming from different lines of research. In line with SSP, we will only consider situations in which both low- and high-frequency information is, in principle, available to the observer, a precondition for filtering. This precludes phenomena such as {\em spontaneous trait inferences} \cite{uleman2008spontaneous} in which observers  “create”  (rather than select) abstract representations from concrete, high-frequency information (for a general theoretical analysis see \cite{gilead2020above}). 

\subsubsection{High vs. low-pass filters}\label{high vs. low}

The leading theory dealing explicitly with filtering of stimulus ranges is {\em Construal Level Theory} (CLT; \cite{trope2010construal, trope2011construal}). CLT posits that people either focus on details or on the whole, a difference often illustrated by the tree vs. forest metaphor. Under low-level construal, people focus on the concrete and specific details of a situation, whereas under high-level construal the broad, abstract, and decontextualized features of the same situation come to the forth. For instance, when planning an event, one may either focus on the nitty-gritty of implementation ({\em how}) or on its larger significance ({\em why}). In signal processing terms, these differences between concrete low-level and abstract high-level construal can be described as applying a high-pass vs. low-pass filter. Thus, although both specific and broad features are theoretically available, the observer selects a specific range on which to focus. How this selection occurs and how it is related to factors such as psychological distance will be discussed in  Subsection IV\ref{third and fourth}. 

The selective filtering of high- vs. low-frequency harmonics has another interesting implication that appears compatible with CLT. Low-pass filters enhance or leave unaltered  low-frequency harmonics while attenuating or deleting high-frequency ones; high-pass filters do the opposite. As low-frequency harmonics correspond to general/abstract categories and high-frequency harmonics correspond to individual characteristics, a person using a low-pass filter will notice fewer interindividual (and possibly more intergroup)  differences than one using a high-pass filter. As a logical consequence, people using high-pass filters should focus on how stimuli are different from each other, whereas those using low-pass filters should focus on similarities. This aligns with claims that, under high-level construal, people “tend to focus less on the multitude of incidental features that distinguish one object from another, and instead focus on the few essential or important features that tie them together”   \cite[p.677]{henderson2013seeing}. In support of this prediction, people under high construal (low-pass filter) have been found to classify objects in broader categories \cite{liberman2002effect} and are more likely to include atypical exemplars into a category \cite[Study 2]{wakslak2006seeing}, presumably because of their greater focus on similarity. 

The filtering of specific frequency ranges also emerges during impression formation as evidenced, in particular, by the literature on {\em individuation vs. social categorization} and stereotyping (for an overview see \cite{bodenhausen2012social}). Researchers investigating first impressions have repeatedly emphasized that observers tend to focus either on individualizing features or on broader categories to which targets belong. For instance,  Fiske and Neuberg (1990, p.2),  proposed a continuum of impression formation from “individuating processes that use a target’s particular attributes to the relative exclusion of category membership” to “category-based processes that use a target’s particular attributes to the relative exclusion of category membership”. This looks much like the high- and the low-pass filters of the upper portion of Figure 3. The choice of a high-pass or a low-pass filter has well-established social-cognitive consequences. Low-pass filters, in which the category membership is the main focus, generally lead to stereotyping, as already available information about the category becomes highly accessible. Similar to what we have seen for CLT, low-pass filters (categorization) also tend to create a focus on similarity rather than difference \cite{park1982perception}, whereas high-pass filters (individuation) promotes a focus  “on how the target person differs from other persons, rather than on class equivalencies within a given group” \cite[p.325]{bodenhausen2012social}, \cite{forster2008effect}. 

A similar distinction between individuation (high-pass filter) and categorization (low-pass filter) can also be found in {\em face processing} and, in particular, in the so-called own race effect. One of the most robust findings in social cognition is that people have poorer memory for faces of outgroup members. These include such diverse categories as race, age, nationality, and even artificially created minimal groups (for a review see \cite{meissner2001thirty}). Hugenberg and collaborators explain this bias in their Categorization-Individuation Model arguing that faces are processed either through individuation or categorization \cite{hugenberg2010categorization}. Individuation means that observers distinguish single faces within a given category, whereas categorization implies that observers focus on the dimensions that define category membership (e.g., race), which naturally implies that they are less able to distinguish one individual from the other during subsequent recognition. We will come back to this model in Section \ref{what filter and why} when discussing the conditions that favour different types of filters.

The fourth broad area of research in which filtering is apparent is {\em language}, whose role  in social cognition is widely known. When describing events, sharing impressions, telling stories or voicing opinions speakers select, from a vast vocabulary, those terms that they deem most appropriate for their current communicative goals. One dimension on which words vary systematically is abstraction \cite{coltheart1981mrc,della2010beyond}. Compared to concrete words (such as table, dog), abstract words (such as democracy, justice) are acquired at a later age, learned mainly through verbal stimulation rather than direct experience with the object, and are more difficult to envisage or contextualize. But abstraction not only varies (semantically) within linguistic categories (e.g., nouns), but also between categories. Semin and Fiedler’s (1991) {\em Linguistic Category Model} distinguishes four large classes of interpersonal predicates varying from the most concrete, namely descriptive action verbs such as “to kick”, over  interpretative action verbs such as “to hurt”, over state verbs such as “to detest”, up to the most abstract form, namely adjectives such as “aggressive” (\cite{semin1991linguistic}; for an extension to category nouns such as “athlete” or “Jew”, see \cite{carnaghi2008nomina}). Although all four terms of this example would accurately describe an episode of interpersonal aggression, the psychological implications vary greatly across linguistic classes. For instance, the more abstract the description, the more it implies stability over time and across targets; compared to a person who “kicked” somebody, an “aggressive” person is more likely to behave in similar ways also towards others and in other situations. In our SSP framework, we can consider the concrete end of the language continuum as resulting from a high-pass filter and the abstract end from a low-pass filter. By choosing abstract predicates, speakers restrain their communications to the low-frequency range, which in turn will affect the impressions listeners are likely to form of the episode. Thus, by focusing on the low-frequency range, while filtering out the specifics of the high-frequency range, speakers create an impression of the episode as broad, general, typical, and enduring (for a similar model, concerning action identification, see  Subsubsection IV\ref{default}).

	All four research areas described above (CLT, individuation vs. categorization in impression formation, own race bias, linguistic category model) share the idea that people process social signals in a highly selective way by applying either a high- or low-pass filter while filtering out the remaining frequency range.  In addition, filtering is also incorporated in a number of classical experimental tasks typically employed in (social) cognition research, such as the card sorting task, the “who-said-what” paradigm, chunking, and the {\em Navon Task}. In the {\em Free Sorting Task}, participants receive a number of cards containing specific pieces of information and are asked to sort them freely into as many piles as they see use \cite{blanchard2016evidence}. This task provides not only information about the latent dimensions guiding categorization such as sex and age, see \cite{brewer1989primacy}, but also about the degree of abstraction (frequency band) that people prefer, given that they can organize the stimuli at any frequency level (from a separate pile for each stimulus card - i.e., the highest frequency -  to a single pile under which all stimuli are subsumed - i.e., the lowest frequency). 
	
Similarly, the {\em Who-Said-What Paradigm} \cite{taylor1978categorical} was designed to unobtrusively assess spontaneous categorization. Here, participants receive information about statements made, for example, by three men and three women, or by three Black and three White individuals, on a topic unrelated to gender or race. Subsequently, they are shown a list of all statements and asked to assign each statement to the correct individual. The fact that people generally make many more within category than between-category errors suggests that they inadvertently organize information around gender or race categories, thus employing a {\em low-pass filter}, even where these categories are of no use to the task at hand. Both tasks, free sorting and who-said-what, show the tendency to (re-)organize high-frequency input at a lower frequency level, suggesting the application of a low-pass filter. 

The same principle can be found in {\em Chunking Tasks} in which participants create higher-order representations by combining smaller segments into larger meaningful units \cite{newtson1973attribution}. For instance, in the case of chess, beginners focus on the moves of single chess pieces, whereas expert players perceive and process larger meaningful patterns, which provides a remarkable memory advantage (for a summary of chunking research see \cite{gobet2001chunking}). 

A final example of experimental tasks representing differential filtering is the {\em Navon Task}, also used as a manipulation of construal level. The {\em Navon Task} induces global or local processing by exposing people to composite shapes, such as large letters composed of small letters. High- vs. low-pass filtering can be induced by having participants read out either the small or the large letters on repeated trials. Together, these examples illustrate that different experimental tasks used in social-cognitive research can be interpreted as assessing or inducing low- vs. high-pass filtering. 

\subsubsection{Beyond low- vs. high-pass filters}\label{beyond}

So far, we have discussed high- vs. low-pass filters that either enhance low-frequency harmonics while attenuating/deleting high-frequency harmonics, or vice versa (see upper portion of Figure 3). However, as can be seen in the lower portion of the same figure, there are other types of filters that deserve attention. On the one side, the {\em bandpass filter} preserves a certain range of frequencies, while eliminating frequencies below and above the critical frequency range. Our visual system, when we stand very close to an image, can be represented as a band-pass filter, see Subsection IV\ref{analogy} below. An example of this type of filter in  cognition is Eleanor Rosch’s {\em Basic Level Category Model} according to which natural objects are preferentially categorized at a basic and primary level such as “dog” rather than “mammal” or “border collie” \cite{rosch1973natural,rosch1976basic}. Basic level categories (“dog”) correspond to terms that children learn first and that are used most frequently in any given language. Importantly, these categories are more informative than both, superordinate categories (“mammal”) that contain only few attributes and subordinate categories (“border collie”) that represent too many attributes to be of practical use.\\ 
On the other side, Figure 3 (right lower portion) also represents a type of filter in which both high- and low-frequency ranges are maintained, while the central frequency range is blocked ({\em band-stop filter}). This corresponds to a third possibility identified by Fiske and Neuberg  \cite{fiske1990continuum} in addition to low-pass and high-pass filters. Here, “people also make sense of other people by combining category membership and target attributes to recategorize others; in such cases, impressions are intermediate between fully category-based and fully individuating impressions” \cite[p.1]{fiske1990continuum}, which looks much like the band-stop filter in Figure 3. In this case, people may process individuating information in the light of a person’s category  memberships. An example is a classical study by \cite{kunda1993stereotypes},  in which observers interpreted ambiguous individuating information in category-congruent fashion, suggesting that they used low-frequency information to interpret high-frequency information. Similarly, Cohen \cite{cohen1981person} found that recall for details shown in a video of a woman, described either as waitress or as librarian, aligned with her category membership. For instance, people remembered her having a beer when the woman was described as a waitress, but as wearing glasses when described as a librarian. Thus, both the interpretation and the recall of specific (high-frequency) information is channeled by low-frequency category information, attesting to the simultaneous use of low and high-frequency information. 

A particularly interesting type of band-stop filter are {\em notch filters}, whose defining feature are their particularly narrow stopbands. They are, for instance, used in instrument amplifiers, especially for acoustic instruments such as acoustic guitar, mandolin, etc., to reduce or prevent audio feedback (Larsen effect), while not appreciably changing the rest of the spectrum. In social cognition, thought suppression may be envisaged as a notch filter that blocks out a specific narrow frequency range not to be considered during cognitive processing, such as the “white bear” in Wegner’s groundbreaking research \cite{wenzlaff2000thought}.  There, the effort to remove thoughts about a trivial topic (white bear) is later ironically followed by an increase of such thoughts. Given the frequently found rebound effects, 
the very attempt to temporarily inhibit thoughts about a certain object may be interpreted as short-lived notch filtering.


Although we will focus mainly on high and low-pass filters in this paper, the above examples illustrate that the multitude of filters of SSP can accommodate more diverse social-cognitive phenomena than this review suggests.

\section{What Filter and why}\label{what filter and why}

\subsection{Analogy with the human visual system}\label{analogy}
There is a considerable analogy between this new paradigm for social cognition and our visual system \cite{schyns1998development,schyns1999dr,pylyshyn1999vision,serre2007feedforward,meese2009tutorial,doutsi2018retina}. As a case in point of this analogy, we focus here on the processing and encoding of spatial information in the two-dimensional retinal image (image data arriving ‘bottom-up’ from  our sensory apparatus). In a first coarse description, our visual system is often depicted as an isotropic low-pass filter with a cutoff frequency\footnote{The cutoff frequency  is defined by the property that the energy of harmonics with  higher frequency  is considerably attenuated. This attenuation can be made precise in specific classes such as electronic systems ($-3dB$ of the nominal passband value).}  that decreases with distance \cite[p.525]{oppenheim1997signals}. Indeed, a square painting can be represented by a double sum (double Fourier series, see Section \ref{what filter and why}) such as
\begin{equation}\label{2DFS}z(s,t)=\sum_{j,k} a_{j,k} x_j(s)y_k(t), \quad t_1\le s,t \le t_2.
\end{equation}
Again some of the canonical signals $x_j$ and $y_k$ in the representation have low-, others high-frequency. An isotropic low-pass filter attenuates high-frequency components in the same way in all spatial directions. The cutoff frequency of the visual system increases as we get closer to the image. Indeed, if we stand very close to a painting, all we see are the patches/spots left by the brush strokes such as two-dimensional coarse pixels (very high-frequency = very rapid changes). In such situation, \cite{meese2009tutorial,doutsi2018retina}, the filter resembles more a {\em band-pass filter}, see the third picture in Figure \ref{IF} for the Magnitude Response and  Subsubsection III\ref{beyond}. Such a filter cuts the low-frequency harmonics, here smoothed shapes, which we do not see standing close to the painting. Thus, a better description of our visual system is that it acts in some conditions like a low-pass filter and in others as a band-pass filter as already observed in \cite{shapley1978effect}\footnote{Notice, however, that our visual system (as our hearing system) cannot process harmonics beyond a certain frequency. Thus, when the passband is in high-frequency, the band-pass filter acts similar to a high-pass filter.}. 

As we start moving away from the painting, the passband of the filter shifts to the left allowing us to notice the contour of object/persons (high-frequency = rapid changes) . Moving further away, the band-pass filter becomes a low-pass filter and we finally see what the painter meant us to see: smoothed shapes such as houses, persons, etc. in a coherent  frame (low-frequency = smooth changes). At this final stage, we no longer see the fine texture of the painting as the high-frequency harmonics are deleted by the low-pass filter\footnote{When we stand distant form the painting, our visual system, resembling a low-pass filter, acts like an {\em anti-aliasing filter}, namely a filter which removes the high-frequency content of a signal to be sampled and then reconstructed from the samples. Aliasing is a degenerative phenomenon in the reconstruction due to undersampling, cf.  A.3. }. As an example, a few years back, a {\em hybrid image} was produced by Aude Oliva at MIT for the March 31st 2007 issue of New Scientist magazine combining a Marylin Monroe picture in low-frequency with an Albert Einstein picture in high-frequency, see e.g. \cite[Fig.1]{oliva2017hybrid}. 
A video-clip is also available\footnote{https://www.youtube.com/watch?v=gfvMU36fgKw}. Thus, from a distance, one  sees Marylin Monroe and, getting closer, Albert Einstein. Hybrid images have been used in visual perception for some time, see e.g. \cite{schyns1999dr}.
While signal processing has long been applied to the human visual system, see also, \cite{blake2005role, braje1995human,meese2009tutorial,doutsi2018retina}, we are not aware of any attempt to extend it to social-cognitive processes. 

\noindent
{\em Summary:  Our visual system  may be described as a low-pass filter or a band-pass filter depending on whether the distance from the object is large or small, respectively. }

\subsection{The third  and fourth hypotheses: Passive vs active filtering and social variables affecting the choice of the filter}\label{third and fourth}
{\em We advance a further fundamental hypothesis: Human beings  use {\em passive} low-pass filters, like the oscillator associated to} (\ref{Newton}),  {\em and {\em active} high-pass filters, which require an ancillary supply of energy. 
Overall, people will prioritize low-pass filters because of their greater energy efficiency, especially when social distance is high. As social distance increases, the lower cutoff frequency of the band pass filter decreases} ({\em it is inversely proportional to the social distance}) {\em leading eventually to a situation where the filter resembles a low-pass and only low-frequency harmonics are relevant to the observer’s eye  or mind. The reverse process occurs when social distance is decreased. In this case, observers will shift the filter towards the high-frequency range, but only if sufficient motivation and cognitive resources are available to provide the necessary ancillary energy.}

In Section \ref{representation}, we have argued that, although generally known under different names, the filtering of information in a specific frequency range is a common and well-documented process in social cognition. To prevent information overload, people (often referred to as {\em cognitive misers}) typically focus on a specific range of stimuli while filtering out stimuli in frequency ranges that are not relevant to their current aims or motivational state. But what motivates people to select a low- or high-pass filter? We argue that, everything else being equal, low-pass filters are the default option as they save energy. Different variables, however, may shift filtering towards the high end of the continuum, including, but not limited to, four relevant variables, namely distance, motivation, cognitive resources, and culture, which will be addressed below.  


\subsubsection{Low-pass filter as default}\label{default}

   As we have seen in  Subsection III\ref{frequency selective}, low-pass filters associated to damped harmonic oscillators (\ref{Newton}) have a number of advantages: they have better stability properties, do not require a supplementary power supply,  and  reduce the risk of errors in reconstruction  due to under-sampling, see the aliasing effect discussed in Subsection IV\ref{analogy}. If these low-pass filters require less samples, are less demanding and more stable than high-pass filters, then, anything else being equal, observers should prefer low- over high-pass filters. In line with this idea,  social psychologists have argued early on that categorization (hence, low-pass filters) often constitutes the default option\footnote{Here we are not treating the issue of which attribute takes on the role of superordinate category and why, a question dealt with elsewhere \cite{bodenhausen2012social}.}. Already Allport (1954) had proposed that categorization has priority due to the need to simplify complex social environments. Similarly, Fiske and Neuberg (1989, p.84) posit that “people attempt processes at the category-based end of the continuum (and are frequently successful) before attempting more individuating impression formation processes. In that sense, category-based processes have priority over more individuating processes” \cite{fiske1989category} . The reason resides in their greater efficiency, given that “category-based processes minimally deplete cognitive resources, compared to attribute-oriented processes” \cite[p.15]{fiske1990continuum}. This initial idea has stood the test of time, considering that recent theorizing is largely in line with earlier models. For instance,  Macrae and Bodenhausen (2000) posit that  “rather than considering individuals in terms of their unique constellations of attributes and proclivities, perceivers prefer instead to construe them on the basis of the social categories (e.g., race, gender, age) to which they belong” \cite[p. 95]{macrae2000social}. Given the intrinsic link between categorization and stereotyping, it is not surprising that a similar argument has been made for stereotyping. For instance,  Neuberg and collaborators \cite[p.250]{neuberg2020toward} state that “one can reasonably argue that stereotyping is the default process through which people come to understand one another, at least early in acquaintanceship”.  Finally, the same principle seems to apply to the way people conceptualize actions. Similarly, Vallacher and Wegner \cite{vallacher1987people}, the authors of action identification theory, argue that “when both a lower and a higher level act identity is available, there is a tendency for the higher-level identity to become prepotent”, implying that people “are always sensitive to the larger meanings, effects, and implications of what they are doing” (p. 5). Translated into SSP terms, we could say that, when both high and low frequencies are available, low-pass filters will prevail. As these examples illustrate, there seems to be broad agreement across different social-cognitive theories that low-pass filters are generally the default option or the starting point, although some authors have challenged this assumption in recent years \cite{monroe2018re}. Once more, this is in line with the  {\em passive} nature of oscillators-related low-pass filters, see the Summary at the end of  Subsection III\ref{frequency selective}. This preference for low- over high-frequency should increase as resources become more limited, as we will see under Subsubsection IV\ref{motivation}. There is still some disagreement whether initial categorization constitutes a truly automatic process (see \cite{macrae2000social} for a discussion), however, most authors seem to agree that categorization requires fewer cognitive resources than the processing of individuating information. This is in line with the idea that high-pass filters are active and require   supplementary energy above and beyond the energy needed for low-pass filters.
   
Another reason why low-pass filters of the form considered in  Subsection III\ref{frequency selective} may constitute the default option is their better stability properties with respect to a high-pass such as the differentiator also described there. Recall that BIBO-stability for the family of filters we consider can be described as the property that small perturbations of the input cause small perturbations of the output. We proceed next to show that this  property has also been observed in the social cognition literature. Indeed, initial impressions of both, groups and individuals, are formed rapidly, but often subsequently encountered new information does not perfectly corroborate the initial impression. In line with the SSP interpretation, such new information should affect the impressions of individuals (high-pass filter) more than those of groups (low-pass filter) \cite{hamilton1996perceiving}.  Hamilton and Sherman explain this pattern by the fact that people expect and perceive greater unity, coherence, and consistency for individuals than for groups; as a consequence, new information not matching the first impression is more puzzling and has a greater weight in this case. In line with this reasoning, there is considerable evidence that the same inconsistent information elicits greater attention, a more intense search for causal explanations, and better memory when referring to individuals than to groups or categories. As a logical consequence, updating of the initial impression in face of new evidence is, on average, more likely for individual targets, whereas impressions of groups or categories remain more stable \cite{wyer1984cognitive}. This is not to say that category representations are, in principle, stable. In line with recent connectionist models, there is increasing evidence that both individual and category representations are highly sensitive to context, resulting in considerable variation over time (for a discussion see \cite{garcia2006stereotypes}). We are simply arguing that there is a relative difference in stability, with category or group impressions, obtained through a low-pass filter, being more stable than impressions of individuals or exemplars, obtained through a high-pass such as a differentiator.

At the same time, social-cognitive psychologists have also argued that, although low-pass filters are often the starting point, there are specific conditions that motivate or force observers to abandon low-pass filters. First of all, observers may shift towards higher-pass filters as familiarity with the stimulus person increases. For instance, Taylor \cite{taylor1981categorization} argued that “as one becomes more familiar with a given group, categories of subtypes will develop. This process will, in turn, lead to a more highly differentiated set of stereotypes for that social group” (p. 87). It is worthwhile to stress here the striking analogy with recent studies of the Outer Plexiform Layer of the human retina \cite[p. 3484]{doutsi2018retina}. According to the authors, the retina-inspired filter: 
"...when it is applied to a still
image {\em visible for a given time}... varies from a {\em lowpass filter}
to a {\em bandpass filter}".  
The latter, as observed in Subsection IV\ref{analogy}, includes higher frequencies.


A similar shift towards the high-frequency end occurs when people encounter information that challenges the low-frequency information. According to some authors, early categorization provides the very basis for noting deviations from category-based expectancies \cite{macrae2000social,sherman1998stereotype}. Facts that are incongruent with existing category knowledge often lead observers to shift attention towards higher frequencies, as in the case of subtyping. When observers encounter others who do not match the stereotypes associated with their category membership, they tend to engage in subtyping or subgrouping (for an overview see \cite{richards2001subtyping}), thereby increasing the frequency range of the filter.

The idea that people shift the filter from the low- to the high-frequency range is also at the heart of action identification theory \cite{vallacher1987people,vallacher2014theory}. According to this theory, any action (e.g., locking the door) can be construed either in terms of the details of the action (e.g., putting a key in the lock) or in terms of its broader meaning, including the actor’s motives  and goals (e.g., securing the house) \cite{vallacher1989levels}. Although people are likely to initially prefer to construe actions in terms of their high-level identities, these are abandoned in favor of lower levels of identification “in the face of high-level disruption” (p. 5). Translated into SSP terms, low-pass filters are preferred unless obstacles are encountered that force the observer to focus on higher frequencies. Together, these different bodies of literature seem to agree that observers often apply low-pass filters as the first and preferred option, but they shift filters in the direction of high frequencies as familiarity with the stimulus increases or when encountering obstacles or contradictions that cannot be resolved in the low-pass mode. 

We shall now examine four factors, extensively investigated in social-cognitive research, that determine the nature 
of the filter, i.e. its passband. These are psychological distance, motivation, cognitive resources and culture. Since motivation and cognitive resources are, by definition, intertwined, they will be discussed under a common heading.

\subsubsection{Distance}\label{distance}
The first variable that regulates filtering is (psychological) distance. This argument is at the heart of CLT \cite{trope2010construal, trope2011construal}. CLT is a general theory of psychological distance arguing that the farther removed an object is from the here and now, the more abstract the mental construal of that object. Vice versa, high-level construals will evoke more distant objects. Abstract mental construals are defined as relatively abstract representations in which central, global features of the object are prevalent. Abstract construal is typically assessed by measures such as language abstraction, breadth of (superordinate vs. subordinate) categories created in verbal or sorting tasks, reference to large vs. small   units of analysis, attention to  global vs. local  features and the like (for an overview see \cite{trope2010construal}). As distance increases, construal becomes more abstract and this holds for perception, categorization, action identification and person perception alike \cite{trope2010construal}. 

Resembling filtering in the visual realm (see Subsection IV\ref{analogy}), CLT  contends “that people use increasingly higher levels of construal to represent an object as the psychological distance from the object increases.” \cite[p.441]{trope2010construal}. Remarkably, this holds for different, strictly interrelated distance dimensions \cite{bar2006association,fiedler2012relations}, including spatial, temporal and social distance, as well as hypotheticality (i.e., whether an event is likely or unlikely). As any of these distances increases, construal becomes more abstract or, in SSP terms, low-pass filters prevail. To cite only a few empirical examples, NYU students described the same event in more abstract language when it was taking place on the NYU campus in Florence, Italy, rather than on the NYU campus in New York \cite{fujita2006spatial}. Similarly, objects were grouped into larger categories when they served for events that were going to take place in the distant rather than near future \cite{liberman2002effect}. Analogous effects were found for hypotheticality or event likelihood. When imagining to plan an event (e.g., going on a camping trip), people grouped relevant items (e.g., tent, matches, rope etc.) into broader categories when the event was described as unlikely rather than likely \cite{wakslak2006seeing}. Turning to social distance, which is the most relevant dimension for social cognition, the behaviors of dissimilar others were interpreted mainly on the basis of primary, central features, whereas more weight was given to subordinate, secondary features when judging the same behaviors performed by similar others \cite{liviatan2008interpersonal}. Thus, shifts towards higher-level construals have consistently been found as any of these distances (space, time, hypotheticality, social distance) increases. 

The parallelism between the visual and the cognitive realm is striking. Looking from the distance at a landscape, such as a forest or a skyline, observers are unable to distinguish single constituents. But when getting closer, they can differentiate trees or buildings, each with its own peculiarity. Similarly, analyzing social environments at high levels of construal, one notes homogenous groups of people, as in the case of the outgroup homogeneity effect \cite{boldry2007measuring,ostrom1992out}. Increased social closeness is needed to appreciate the differences between single members of a social group. An emblematic case of such social closeness is ingroup (vs. outgroup) membership. Generally, people perceive greater differences between members of their own group than between members of outgroups. Whereas the ingroup is perceived as variegated and complex, the outgroup appears undifferentiated and homogenous. Although different explanations have been offered for this {\em outgroup homogeneity effect} (including differential familiarity), a prominent approach contends that ingroup and outgroup information is processed in qualitatively different ways, namely through greater individuation (differentiating unique individuals) or subtyping in the case of ingroups \cite{park1982perception}. In line with this general idea, compared to outgroup labels, ingroup labels activate more subcategories \cite{park1982perception}. According to this differential processing hypothesis, the same information is processed in qualitatively different ways, which, in our opinion, greatly resembles low- and high-pass filters. 

An analogous difference has been shown for {\em self} vs. {\em others}. For instance, people tend to describe their own behavior as "depending on the situation" but readily ascribe global traits to others, cf. \cite{goldberg1981unconfounding,jones1987actor}. Self (compared to other) construals tend to be both, more concrete \cite{bruk2018beautiful} and more complex as they also contain non-observable information derived from introspection such as hopes, fears, and emotions \cite{karylowski2020temporal}. This self-other distinction is not categorical, but varies dynamically across situations. For instance, the difference in perspectives on self vs. others diminishes as time distance increases (e.g., considering the self in the distant past) and as familiarity with others increases \cite{karylowski2020temporal}. In signal processing terms, the {\em self} is, on average, filtered through a higher pass filter than {\em others}, but this difference between self and others varies in theoretically meaningful ways. As argued by some authors \cite{park1990measures,simon1993asymmetry}, the individuating perception of the self may, at least in part, drive the outgroup homogeneity effect, given that ingroups, but not outgroups, contain the self as a unique element. Thus, the two phenomena not only rely on similar processes (high- vs. low-pass filtering), but one may also contribute directly to the other.
	Together, the above lines of research, CLT, outgroup homogeneity effect, self vs. other construal, support the idea that closeness is, in SSP terms, associated with high, and (social) distance with low-pass filtering. We will come back to this issue in  Subsection V\ref{comparison}, where we will systematically compare CLT and SSP.
	
\subsubsection{Motivation and cognitive resources}\label{motivation}

    Processing high frequency information is taxing and requires engagement with the task. Social cognition approaches have acknowledged early on that piecemeal integration of individuating bits of information require greater motivation and cognitive resources than category-based information processing. We have seen in  Subsection II\ref{first hypothesis}, that negative information is, in SSP terms, of higher frequency and positive information of lower frequency \cite{alves2017good}. Although low-pass filters are generally the default, this is not the case here. There is ample evidence that people pay greater attention to negative information and process it more thoroughly than positive information. This asymmetry in favor of the negative information is generally accounted for by the fact that negatively valenced stimuli are of particular relevance \cite{baumeister2001bad}. From an evolutionary perspective, it is adaptive to pay attention to negative stimuli, especially when they are dangerous or threatening. As negative events signal a need for change, and often a need for action, it is not surprising that people generally apply high-pass filters, attending more to negative than to positive information. The reasoning here could also be reversed, suggesting a feedback loop. Specifically, we may argue that the motivation triggered by dangerous or threatening situations prompts high-pass filtering, and that negative information per se carries high-frequency components (see Subsection II\ref{first hypothesis}), therefore reinforcing high-pass filtering.  
    
Another prototypical example of how motivation and personal relevance affects filtering is Petty and Cacioppo’s (1986) {\em Elaboration Likelihood Model}. This model posits that  a persuasive message may be processed through one of two routes: Under the central route, people carefully scrutinize the quality of the arguments presented in favor of or against a certain position, whereas under the peripheral route people will take a shortcut by relying on “proxies” ({\em peripheral cues}), such as the status, attractiveness or credibility of the source or the length of the message (\cite{petty1986elaboration}, for a recent extension see \cite{kitchen2014elaboration}). Whether people will engage in  central (rather than peripheral) processing depends on their motivation and the availability of cognitive resources. For instance, people will be more likely to choose the central route when the topic is of personal relevance (thereby increasing motivation) and when sufficient time and energy are available to process the arguments. Motivation and cognitive resources operate in a multiplicative fashion, such that central processing becomes unlikely whenever one of the two factors approaches zero. In other words, for the central route, people have to be both, motivated and able to process the message. The model predicts a tradeoff between a superficial but speedy peripheral route and a more accurate but more demanding central route, which may be equated with low- and high-frequency filters, respectively. A similar argument can be made for Chaiken’s  {\em Heuristic Systematic Model} \cite{chaiken1987heuristic}, which closely resembles Petty and Cacioppo’s elaboration likelihood model \cite{petty1986elaboration}, but with an important difference: the two elaboration modalities can co-occur. In SSP terms, whereas the filter in Petty and Cacioppo’s model shifts towards the high frequency range as motivation increases, in Chaiken’s model high and low frequencies may, at times, be processed simultaneously, resembling a band-stop filter (see  Subsubsection III\ref{beyond} above). More recently, such dual process models, have been challenged by unified theoretical approaches such as the {\em Unimodel} \cite{kruglanski1999persuasion,kruglanski2011intuitive}. The unimodal posits that there is no inherent difference between peripheral cues and message arguments in the persuasion process. Both peripheral and central information contained in messages is processed to the degree that motivation, information relevance and processing capacity are high. In other words, how much (cue- or message-related) information the person is willing and able to process depends on variables such as self-relevance and motivation. In SSP terms, this resembles the band width of the filter, such that motivation increases the processing of both peripheral, low-frequency cues and message-specific, high-frequency cues. If this interpretation is correct, then SSP seems to offer an overarching perspective able to accommodate all three persuasion models by considering specific filter properties, namely  the passband of the filter ({\em Elaboration Likelihood Model}), the possibility of a band-stop filter ({\em Heuristic Systematic Model}), and the band width of the selected filter ({\em Unimodel}).

Similarly, when developing their influential {\em Impression-Formation Model}, Fiske and Neuberg (1990, p.8) argued that piecemeal integration will take place only when less demanding category-based processes have failed and when people have “sufficient time, resources and motivation” to proceed in an attribute-by-attribute fashion \cite{fiske1990continuum}. Thus, observers will go beyond the initially activated category (e.g., gender) and attend to additional individuating information (e.g., profession, hobbies, possessions, habits) only when they are motivated to do so, for instance because the target person is personally relevant to them, because the interdependence structure is such that their outcomes depend on the target, because there are accountable, or because of numerous other reasons motivating observers to draw accurate inferences (for an overview of such additional processes see \cite{macrae2000social}). Importantly, seemingly identical situations can lead to opposite filtering depending on the activated goal. For instance, competitive interdependence between individuals increases the likelihood of individuation \cite{neuberg1987motivational}, whereas competitive interdependence between groups increases categorization \cite{brewer1995person}.

Analogous effects of motivation are apparent in the processing of same- vs. other-race faces. We have already seen that people are better at recognizing own- than other-race faces \cite{meissner2001thirty,timeo2017race}. This is, at least in part, related to the lesser effort put into individuating other-race faces. For instance, the same faces are processed more holistically and remembered better when arbitrarily defined as ingroup members \cite{hehman2010division,van2012social}, presumably because people allocate greater attention to faces labelled as belonging to the ingroup. This suggests that shared group membership motivates people to attend to individuating information rather than treating faces as interchangeable exemplars of a category. 

However, the effortful processing of social signals through high frequency filters not only requires greater motivation, but also more cognitive resources. Again, research on both, social categorization vs. individuation and on the own-race effect, provides evidence for this assumption. As far as categorization and individuation concerns, the reliance on categories (low-pass filters) rather than on individuating information (high-pass filters) is more likely when people are under time pressure \cite{dijker1996stereotyping}, ego-depleted \cite{govorun2006ego}, or at a nonoptimal time of their circadian rhythm \cite{bodenhausen1990stereotypes}. A particularly telling example of the role of cognitive resources in filtering is a classical study by Rothbart and collaborators \cite{rothbart1979recall}. Here participants received trait information (e.g., lazy, polite) about single individuals (e.g., Joe, Bill) who were part of a fictitious group (Dallonians).  The main result of the study was that people organized the information around the characteristics of individual members when memory load was low, but around the group as a whole when under high memory load. In SSP terms, when cognitive resources were scarce (high memory load), people used a low-pass filter to facilitate their work. Turning to the own-race effect, T\"uttenberg and Wiese analyzed event-related brain potentials and found that the encoding of other-race faces requires additional cognitive resources in order to be processed similar to ingroup faces \cite{tuttenberg2021recognising}.

A final example that may be interpreted under a SSP perspective are moral dilemmas. In the typical research paradigm, participants are asked to decide whether to engage in hypothetical, immoral behaviors in order to prevent an even greater damage, such as killing of one person to save the lives of many (for instance, trolley dilemma). Decisions generally either follow a general deontological principle or inviolable moral rules (which, in SSP terms, may be considered low-frequency) or a situation-specific utilitarian consideration in which the morality of the action depends on its unique outcome (which, in SSP terms, represents a case of high-frequency). In line with our argument that high-pass filters are more demanding, utilitarian judgements have been found to require greater working memory \cite{moore2008shalt}, a more deliberate cognitive style \cite{bartels2008principled,paxton2012reflection}, and more time \cite{suter2011time} compared to deontological judgments. 

There are many situations in which people apply low-pass filters simply because the task demands exceed their cognitive resources. Consider for example the case of resume screening in hiring procedures with disproportionate numbers of job applicants, a phenomenon known as resume overload. This task may turn out to be so taxing and complex that evaluators rely on low-frequency cues such as age, gender or ethnicity, often leading to hiring discrimination \cite{koch2015meta}. Similarly, in the who-said-what paradigm mentioned earlier (see Subsubsection III\ref{high vs. low}), high-frequency harmonics constitute the only useful type of information to accomplish the task, but low-pass filters are activated because high frequency harmonics are simply too many to be processed in detail. Together, the above examples suggest that accessing  high frequency harmonics requires additional energy, here in the form of motivation and cognitive resources, and that in the absence of such additional energy, low-pass filters will prevail (cf. our third hypothesis at the beginning of  3.2).
}

Does the application of low-pass filters necessarily lead to reduced accuracy and poor decision making? We believe that this is not necessarily the case. Low-pass filters offer a faster and easier solution to solve complex problems in situations in which time or cognitive resources are limited. This view is in line with the bounded rationality approach proposed by Nobel price winners Herbert Simon and Reinhard Seleten, who posit that under specific circumstances cognitive processes deviate from perfect rationality to manage the costs that would be associated with the latter. This argument is taken one step further in Gigerenzer’s research program on fast and frugal heuristics, according to which overcoming rationality is not only adaptive to the context, but may even lead to better cognitive outcomes. If the social decision is to be taken fast, or if the high frequencies are too many or too few or of poor quality (e.g., redundant) to lead to sensible conclusions, a high-pass filter is of no advantage \cite{gigerenzer1999simple,gigerenzer2009homo}. For example, the less-is-more effect shows that a limited amount of knowledge can lead to more accurate decisions when relying on heuristics, which  we could equate to a low-pass filter. For example, in deciding which of two American cities (San Diego vs. San Antonio) was larger, German respondents (with less knowledge of American geography) performed better than American students, because they applied a low-pass filter, namely the recognition heuristic. German students simply distinguished the city they did  know from the city they did not know, as famous cities are likely to be larger \cite{gigerenzer1996reasoning}. For the same reason one should not assume a-priori that stereotypes (low-pass filters) are necessarily less accurate than individuating information (high-pass filters) \cite{judd1993definition,jussim2015stereotype}. Thus, reduced cognitive resources are likely to shift filters towards the low-frequency end, but this does not automatically lead to poorer performance.

\subsubsection{Culture}\label{culture}
The propensity to apply high- vs. low-pass filters, or narrow vs. broad passbands, also seems to reflect cultural differences. There is ample evidence that some cultures tend to process complex information analytically, others more holistically \cite{nisbett2001culture,kitayama2019east}. People from Western cultures tend to use an analytic processing style in which they analyze the attributes of a focal object independently of its context. In contrast, East Asians tend to process the focal object and its surrounding holistically, as a whole composed of interrelated elements. For instance, an underwater scene composed of fishes swimming among algae is processed in a piecemeal fashion by American participants who label and recall single elements (fish, algae, etc.). In contrast, East Asian participants create a unifying, holistic representation of focal objects and background and, as a consequence, are very sensitive to any subsequent changes in the background \cite{masuda2001attending}. Similarly, in the above-mentioned {\em Navon Task}, participants in collectivist cultures show a stronger tendency to attend to global (rather than local) features \cite{mckone2010asia,oyserman2011culture}. Even self-construal seems to follow the same pattern;  while people from Western cultures tend to perceive the self as independent from others and from the social context, East Asians construe the self as closely interconnected with others and as embedded in a social context. Thus, cultural norms either stress the individual specificity of each person acting in a social vacuum or as persons embedded in a complex social context.

In SPP terms, the latter differences between analytic and holistic cultures \cite{nisbett2001culture}, namely local vs. global processing and independent vs. interdependent self construal, seem compatible with the idea that Westerners apply relatively high- and East Asians relatively low-pass filters. In contrast, the former example (underwater scene) points to a somewhat different effect of culture on signal processing.  While Westerners process the focal object at the expense of the context, East Asians process both, integrating them into a holistic representation. In SSP terms, this suggests that Westerners process mainly high frequency signals while neglecting the lower frequency, whereas East Asians process both. If this interpretation is correct, this would suggest that the passband is larger for East Asians. Looking at the upper right portion of Figure 3 (high-pass filter), we can envisage that the blue box extends towards the left (lower frequencies) for East Asians. While the filter includes the higher frequencies for both groups, it extends further towards the lower end for East Asians only. Together, the above lines of  research suggest that some cultures are more likely to promote low-pass filters and/or broader passband than others, although the historical, philosophical, social, linguistic and architectural underpinnings of these cultural differences are still debated \cite{caparos2012exposure, nisbett2005influence, varnum2010origin}.

As this brief and selective review of relevant literature shows, psychological distance, motivation, cognitive resources and culture jointly determine whether people apply high- or low-pass filters, and their passbands, when processing social information. The overview of the four factors is not meant to be exhaustive. There are many other variables likely to affect the level of filtering. On the one side, task characteristics that demand attention to higher or lower frequencies are likely to shift filters accordingly, as shown in both vision \cite{schyns1999dr} and social-cognitive research \cite{fiske1987category}. On the other side, individual differences such as need for closure \cite{dijksterhuis1996motivated}, miserly thinking \cite{blanchar2020individual},  limited executive control \cite{payne2005conceptualizing}, and promotion vs. prevention orientation \cite{forster2005global,pham2005promotion} are likely to increase the reliance on low-frequency categories at the expense of individuation.



\section{Conclusions and Outlook}\label{conclusions}

We believe that the SSP approach offers a promising integrative framework and common language able to accommodate a wide range of social-cognitive phenomena. SSP complements, but does not substitute any of the models discussed here, each of which offers a unique and detailed perspective on a specific phenomenon such as other race recognition, impression formation, language use, and the like. SSP helps to identify some of the common underlying principles of these models and offers a formal language to describe these communalities, without pretending to deal with the specifics of each of these models. SSP describes how frequency ranges are selected. In its current form, it does not deal with the question of how abstract concepts come about in the first place \cite{gilead2020above}, how they come to be represented  in language \cite{borghi2019words}, how they are acquired by children \cite{timeo2017race} or how they are represented in the brain \cite{gilead2014mind}. At the current stage, the predictions of SSP are limited to adult participants, to situations in which high and low frequency information is simultaneously available, together with linguistic labels able to encode the event at different levels of abstraction. We will now discuss similarities and differences with respect to CLT, to then propose new, testable hypotheses, uniquely derived from SSP.

\subsection{Comparison between CLT and SSP}\label{comparison}

\begin{figure}
\begin{center}
	\includegraphics[scale=0.35]{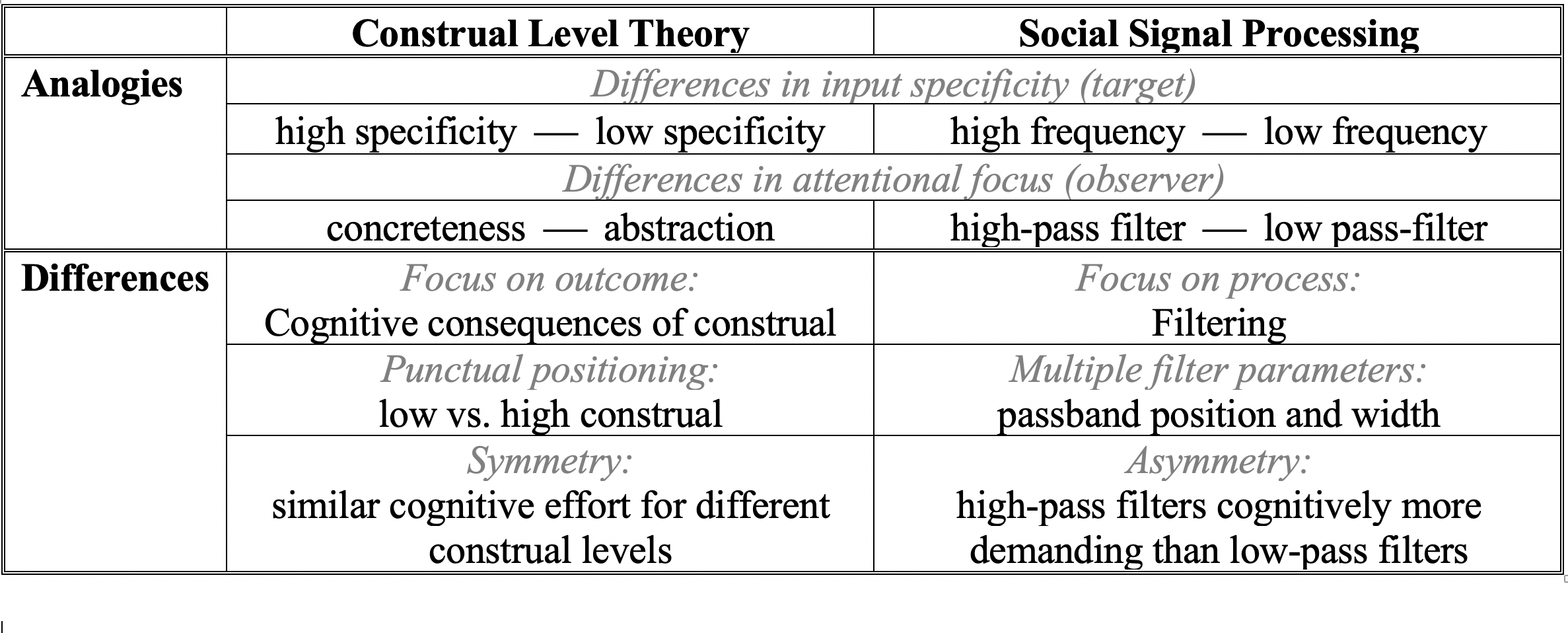}
	\caption{Common and distinct features of CLT and SSP}
	\label{TCR}
	\end{center}
\end{figure}

As is evident from Sections 1 through 3, our SSP approach has multiple points of contact with CLT. It is therefore worthwhile to detail what features the two models have in common and where the differences between CLT and SSP lie. In Table \ref{TCR}, we schematically represent both the common and the distinct features of CLT and SSP.
Two common features have been extensively discussed in Sections I and II, but it may be useful
to recap here their core characteristics. First, both theories position input signals on a continuum according to their degree of specificity. At the maximum level of specificity, we find fine-grained elements, namely {\em the trees} in the classic metaphor used in the CLT framework, corresponding to high frequency harmonics in SSP. On the opposite pole of the continuum, we can envisage {\em the forest}, which, for both theories, represents the most abstract level. The second common element pertains to the attentional focus of the observer. According to both theories, since the information is generally rich of fine-grained elements that merge into larger compounds, the person processing the information can shift attention across different specificity levels depending on contextual or personal drives. In CLT terms, the same stimuli may be processed at a high construal level, which corresponds to a low-pass filter in SSP. To the contrary, the high specificity end of the continuum is captured through a high-pass filter in SSP and results in low construal in CLT. Thus, people may process the same stimuli at different degrees of specificity and the two theories show a reasonable overlap in the identification of the conditions that define which specificity level will be in the spotlight of the cognitive activity. In particular, CLT posits that psychological distance from the stimulus is the major driver of high vs. low construal. Objects that are psychologically distant (in terms of time, space, social distance, or hypotheticality) are likely to be construed at higher abstraction levels, whereas minute details are more likely to be considered when processing stimuli that are perceived as psychologically close. The two theories also concur that the choice of high vs. low construal or low- vs. high-pass filters, respectively, also depends on chronic preferences for a given level of specificity (see e.g. \cite{freitas2001abstract, vallacher1989levels}), on personality, and on motivation (see Subsubsection IV\ref{motivation}).  
Although both approaches stress the importance of stimulus specificity and the flexibility in processing stimuli at different levels of specificity, the two theoretical frameworks differ in at least three important ways, namely 1) their focus on outcomes vs. processes, 2) the positioning of mental construal at a specific point of the continuum vs. passband filtering, and 3) the symmetry (vs. asymmetry) of different levels of input specificity.

The first and most general difference is that CLT focuses on the cognitive outcome of stimulus construal,           with psychological distance being the primary prompt of such construal. In the typical CLT paradigm, psychological distance is experimentally manipulated by enhancing spatial or temporal distance or by reducing familiarity and then the outcome of this manipulation is measured. Differently, SSP’s  scope regards the mechanics of the processes at play. The major goal of SSP is to describe the process that links the manipulation to the outcome, and this process consists in the filtering activity that is cognitively enacted while processing social information. Although we have here focused on the general, conceptual aspects of SSP, the approach offers the possibility to mathematically define the specific parameters that lead to different construals or, in SSP terms, that determine different filtering actions. 

The second critical difference between SSP and CLT lies in the fact that the CLT paradigm generally defines construal level as a specific point on the concrete-abstract continuum. Differently, SSP contains several other parameters. Filters are grouped in families depending on the frequencies that get enhanced or attenuated. These families include low- vs. high-pass filters, similar to CLT, but also band-pass and band-stop filters. Moreover, within each of these families, filters differ in the location of their passband and in the passband width\footnote{Typically, several other parameters are needed to fully characterize a filter which is  described by its frequency response (\ref{FR}).}. For instance, filter passbands may be very narrow, in extreme cases selecting a single harmonic, or they may be rather broad, encompassing a large range of harmonics of different frequency. Psychologists are quite familiar with the idea of band width, starting from Easterbrook’s  groundbreaking hypothesis (1959) that, as arousal increases, observers narrow their attentional focus to a smaller range of stimuli. In fact, threatening stimuli not only capture attention, but also narrow the band width of the filter \cite{sorensen2014threat}. The same principle can be found in social cognition. For instance, when encountering a person belonging to a given social category, the observer may mainly perceive relevant category memberships (low-pass filter), or mainly individuating features (high-pass filter). In both cases, the bandwidth may be either very small or relatively large, for instance, including not only category membership information but also some of the individuating features.  

This difference between the two models is  important, because it allows for specific predictions or explanations of social phenomena. For instance, different from SSP (see Subsubection IV\ref{motivation}), CLT does not predict any systematic differences in construal level as a function of culture. In fact, culture is not considered in Trope and Liberman’s fundamental review \cite{trope2010construal} and two meta-analyses revealed no evidence for culture moderating the effect of psychological distance on construal level, nor on any downstream consequences for judgments, decisions, and behavior \cite{soderberg2015effects}. Possibly, the reason why CLT is not specifically addressing cultural differences is because it focuses  on specific points on the specificity spectrum rather than on the band width. Assuming that any culture can see both trees and forests, we suspect that this punctual appraisal of construal levels is not suited to detect cultural differences. However, as we have argued under Subsubsection IV\ref{culture}, cultures may differ in their preference for narrow or broad band filters, as in the case of East Asian participants who, different from their American counterparts, tend to attend to both, the focal object and the broader context (see binding principle). Possibly, cultural differences do not rely so much on the positioning of the filter along the frequency continuum, but rather on its bandwidth. 
A similar argument can be made for the moral dilemmas discussed in Subsubsection IV\ref{motivation}, where people resolve sacrifice dilemmas (e.g., trolley dilemma), either in line with general deontological principles (low frequency) or with  utilitarian arguments (high frequency); under this dichotomic view, women are generally believed to prefer the former, men the latter type of moral reasoning. However, when using more refined measures \cite{friesdorf2015gender}, it becomes clear that women apply a wider passband filter in appraising the dilemma than men do. These two examples suggest that it may indeed be useful to consider band width in addition to the punctual collocation on the abstract-concrete continuum.

The third key difference between the two theories regards the cognitive cost attributed to different specificity  levels. In CLT terms,  “extracting the general meaning and invariant characteristics of objects is not necessarily more or less effortful than fleshing out the minute details of these objects”  \cite[p.15]{trope2010construal}. In contrast, SSP asserts a strong  asymmetry between high- and low-pass filters (see  Subsections III\ref{frequency selective} and IV\ref{third and fourth}). First of all, such asymmetry is already apparent from the fact that the lowest frequency is well defined (zero) corresponding to the constant harmonic, but it is impossible to define a “highest frequency”. In principle, one can conceive increasing frequencies without any limit, leading to wildly oscillating harmonics. The second, deeper, reason for the asymmetry was explained in detail in  Subsection III\ref{frequency selective}: Effective high-pass filters require “supplementary energy”. This led us to formulate our third hypothesis at the beginning of  Subsection IV\ref{third and fourth}. According to this view, high pass filters demand more cognitive resources and greater motivation. As a consequence, low-pass filters are more promptly available and generally used as default, especially when cognitive resources are limited. 
In conclusion, despite the evident similarities, SPP differs from CLT in its focus on process vs. outcome, in its reliance on multiple parameters (among which the location and the width of the passband), and in the assumption that high-pass filters employed in social cognition require more cognitive resources (energy, in SSP terms)  than lower-pass filters.

\subsection{Derivation of new hypotheses and future research.} \label{new hypotheses}
The main aim of the present paper was to lay down the basics of SSP. We believe that this theoretical framework also offers a number of novel perspectives whose implications remain to be explored, some of which we will outline below.

First and foremost, SSP has a multitude of parameters, that allow us to describe diverse social-cognitive phenomena. For a filter, besides passband location and width, these include absolute value (height) and phase of the frequency response in the passband. Actually, the full description  of a filter requires  its {\em frequency response} (7), a complex-valued function changing the power and the phase (Appendix \ref{real}) of each single harmonic according to (8). 
To avoid excessive complexity, in this paper we have focused mainly on two parameters that, we believe, are of particular interest to social psychologists, namely passband location and width. Whereas passband location is contained in other social-cognitive models (among which CLT), to our knowledge passband width has rarely been considered in social cognition, either in itself or in combination with location. The role of  passband width has been recognized in attention since the 50ies (see Easterbrook hypothesis), offering an important explanation of the Yerkes-Dodson law that predicts an inverse U relationship between arousal and performance \cite{faller2019regulation}. It is likely that accuracy in social cognition (e.g., inferences, predictions) depends not only on the location on the specificity spectrum, but also on the band width of the filter, possibly with best performance at intermediate pass band width. 

Second, the hypothesized asymmetry in cognitive resources required for processing low vs. high frequency stimuli remains to be tested explicitly. Although much of the extant literature is in line with this idea (see Subsubsections IV\ref{default} and IV\ref{motivation}), to our knowledge, there is currently no direct empirical evidence to show that social cognitive processes in the high-frequency range undergo greater interference than those in the low-frequency range when cognitive resources are scarce (e.g., under time pressure, under enhanced memory load, during dual tasks).

Third, SSP has a great potential for future development well beyond the basic processes addressed here. For instance, SSP can be extended to include a time perspective, involving both target and observer, a feature that is often missing in social cognition. On the one side, targets may evolve over time, making different characteristics (harmonics) salient at different times. On the other side, filters are likely to evolve over time, for instance shifting from the low- to the high-frequency end, as need arises (as argued in Subsubsection IV\ref{default}). Such shifts over time can be represented within SSP by time-varying systems such as the {\em adaptive filters} of digital signal processing, \cite{gaydecki2004foundations}.

Fourth, SSP predicts that the low-pass filters humans use are more stable than high-pass filters. Small variations in input can produce large variations at output at high, but not at low frequencies. This has interesting and testable implications for social cognition. For instance, everything else being equal, impressions of social categories (including stereotypes) should show greater stability or invariance over time than impressions of individuals. 

A fifth hypothesis derived from SSP concerns sampling and, more specifically, aliasing, see the Appendix \ref{aliasing}. Aliasing predicts a specific relation between sample size and the inferences that people are likely to draw from such samples. The main idea is that, under-sampling will lead people to infer low (rather than high) frequency harmonics, even when samples are theoretically compatible with both, high and low frequency harmonics. Implications of the aliasing principle for social cognition are interesting as it suggests that people will jump to category (rather than individualizing) inferences whenever judgments are based on few exemplars or few behavior samples. The {\em solo} effect discussed in  Subsection II\ref{first hypothesis} may be interpreted as a specific case of this general principle.

Sixth, the model allows for {\em feedback} to take place.  Although many may not be aware of this, we live in a world where automatic control devices based on feedback are everywhere: from the angular velocity control in our washing machine, to the fuel injection system in our car to the autopilot in an aircraft, to name but a few. It took a long time to understand the many implications associated to closing the loop, making feedback one of the deepest concepts in the history of science, cf. the Appendix \ref{bits}.  While feedback can be used to stabilize a system (see  Appendix \ref{low vs. high filters} for an example), it can, among others, also be used to reduce the effect of disturbances and to make the system more robust with respect to parameter variations.  Although it was impossible in this first paper to explore the role of feedback in social cognition, SSP allows for such future development.  Indeed, social-cognitive processes never occur in a social vacuum, generally evolve over time, and often include feedback loops, yet social cognitive research and theorizing rarely contemplate the existence of feedback systems to be operating (for exceptions see self-fulfilling prophecy \cite{madon2011self,snyder1977social}, the attribution of mental states \cite{zaki2011reintegrating}, and  person construal, \cite{freeman2011dynamic}). SSP may help to close the gap between social psychological theorizing and social-cognitive processes which are, almost by definition, feedback loops.

Finally, if proven accurate, our SSP approach may well lay the basis for understanding the neuropsychological underpinnings of social-cognitive processes, much as has happened in vision. Even in the absence of numerical quantification, SSP offers the conceptual mathematical tools to represent social-cognitive processes. Although, as a field, social cognition is currently quite distant from this stage of development, signal processing-based models are already present in other areas such as vision \cite{petras2019coarse} and computer vision \cite{granlund2013signal}. Thus, it is conceivable that social cognition will follow a similar route in the future. 

In conclusion, we believe that SSP not only offers a promising tool for integrating existing knowledge within an overarching framework, but it also paves the way for several novel inquiries into social cognition.

\appendix\label{appendix}
\subsection{Real and complex numbers}\label{real}
Real numbers can be thought as distance along a line. They admit a finite or infinite decimal expansion. They include rational numbers which can be expressed as the ratio of two integers, but also irrational numbers such as $\sqrt{2}$ (the ratio between the length of the diagonal and of the side in a square) and  $\pi$ (the ratio between the length of a circle and its diameter). The set of real numbers $\R$ is totally ordered in that $a<b$ if $b-a$ is a positive number. The absolute value $|a|$ of a real number $a$ is the non negative real number representing its distance from the origin. Thus, for instance, $|3|=3$, $|-\frac{2}{3}|=\frac{2}{3}$.
\begin{figure}\label{complex}
\begin{center}
	\includegraphics[scale=2]{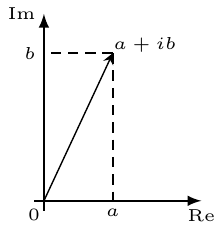}
	\caption{Representation of a complex number in the plane}
	\label{Complex}
	\end{center}
\end{figure}
A complex number can be visually represented as a pair of real numbers $(a,b)$ forming a vector on a diagram  representing the complex plane, cf. Figure \ref{Complex}. ``Re”   is the real axis, ``Im”   is the imaginary axis, and $i$, the imaginary unit\footnote{In electrical engineering, as $i$ is the standard notation for the circuit current, the imaginary unit is denoted by $j$.}, satisfies $i^2=-1$. 

The absolute value $|z|$ of a complex number $z=a+ib$ is again the distance of the tip of the arrow from the origin. By Pythagoras' theorem
\[|z|=\sqrt{a^2+b^2}.
\]
If $b=0$, the complex number is actually real and the two absolute values coincide.

Another important representation of the complex number $z$, alternative to the cartesian coordinates $(a,b)$, is the {\em polar} representation $(|z|,\varphi)$, where $\varphi$ is the angle between $-\frac{\pi}{2}$ ($-180^o$) and $\frac{\pi}{2}$ ($180^o$) formed by the arrow in Figure \ref{Complex} with the $Re$-axis. This representation is often used in signal processing. The distance $r=|z|$ of the tip from the origin $(0,0)$, called modulus or magnitude of $z$, and such an angle, called {\em phase} or {\em argument}  of $z$, determine univocally a point in the plane where the tip of the arrow must be located.

\subsection{Fourier series and filters}\label{fourier}

The unison choir music described in Subsection II\ref{harmonics} we hear  $\{x(t); t_1\le t\le t_2\}$ can be represented  as a {\em weighted sum} ({\em Fourier series})
\begin{equation}\label{FS}x(t)=\sum_k a_k x_k(t), \quad t_1\le t \le t_2.
\end{equation}
Here, the $x_k$ are {\em harmonically related} signals, meaning that they are {\em periodic} of period 
\[T_k=\frac{T}{k}, \quad T=t_2-t_1,
\]
and, therefore, all share period $T$. It follows that their {\em frequency} $f_k=1/T_k=k/T$ is an integer multiple of the {\em fundamental frequency} $1/T$ precisely as it happens with pure tones that are $k$ octaves apart. As for monophonic music, harmonics contributing to the sum (\ref{FS}) differ not only in frequencies but also in {\em amplitude}. Namely, the coefficient $a_k$ of the $k$-th harmonic $x_k$ determines the {\em energy} of that harmonic in the original signal $x$\footnote{More precisely, the $a_k$ are real or complex numbers. Their absolute value squared $|a_k|^2$ is the average power (energy rate) of the $k$-th harmonic of $x$, see Appendix A.1.}. Significant examples of families of canonical signals $\{x_k(t)\}$ are provided next.

In 1807, motivated by his studies on the heat propagation, Fourier advanced the hypothesis that any signal $\{x(t); 0\le t\le T\}$, including those that exhibit a jump in their graph (discontinuous), can be expressed as combination of the form\footnote{The sum may actually be finite with only finitely many of the coefficients $\{b_k\}$ and $\{c_k\}$ different from zero.}
\begin{equation}\label{FSeries}x(t)= a_0 + \sum_{k=1}^\infty \left[b_k\cos\left(2\pi f kt\right) + c_k \sin\left(2\pi f kt\right)\right], \quad f=\frac{1}{T}, \quad 0\le t \le T.
\end{equation} 
Fourier was shown to be correct for a large family of signals ({\em finite-energy signals}) provided, in the case when the sum is actually infinite (series), the convergence of the finite sums is intended in a suitable way (mean-square convergence, see the Riesz-Fisher Theorem e.g. in  \cite{rudin1974real}). This result, and its generalizations, has an important consequence: Basically any signal defined on a bounded domain (finite time interval, a rectangle in the plane) can be represented through a Fourier series. 

Fourier series  have, in the meantime, become a fundamental tool in many areas of science such  mathematical physics and signal processing \cite{oppenheim1997signals,evans1998partial,katznelson2004introduction,sonka2014image}. Other series representations with respect to a family of canonical signals different form sinusoids have also proven to be extremely useful. We mention wavelets series underlying the JPEG 2000 image coding system \cite{taubman2002jpeg2000} as an example.

We close by observing that, as it is apparent starting with Subsection II\ref{first hypothesis}, harmonics appearing in the representation of individuals/social stimuli don't necessarily have frequencies that  are one an integer multiple of the other. For instance, a high frequency characteristics such as sporting blue hair does not lead through the abstraction process (see the language example of Subsection II\ref{harmonics}) to being Black or White. In Signal Processing terms, this means that to properly represent an individual or, more generally, a social stimulus, we need to resort to the {\em Fourier Transform}, see (\ref{FR}) in the following subsection. The latter generalizes the Fourier series in that a continuum of harmonics (one for each point of the real line rather that just one for each integer number) is used to represent the signal.  Moreover, a proper representation of individuals would also involve {\em digital} signals. For instance, in a by now dated classification, sex takes only the two values Female and Male, it is namely a {\em binary} (logic) signal. In conclusion, a more precise description of social stimuli should involve {\em hybrid} signals as sum of both analog and digital harmonics. To avoid obscuring ideas with technicalities, we have left these issues out in this first paper using only  Fourier series. We also mention, as a partial justification of this simplification, that filtering for signals represented through Fourier Series or Fourier Transforms is analogous.

\subsubsection{Convolutional Systems}\label{convolutional}
Consider, for instance, the following problem: We wish to send a voice message to a friend using our cell phone or our computer. How can this be accomplished? Cell towers or internet cannot transmit rapid variations of acoustic pressure (voice). They  can, however, transmit sequences of bits. The first step is carried out by a microphone which transforms the input audio signal $\{x(t);t\in T\}$ into an output (analog\footnote{An analog signal is a continuous representation over time of a quantity, e.g. voltage, which varies continuously.}) electrical signal $\{y(t); t\in T\}$ preserving the waveform (transducer). Here $T$ is the time of the recording, see the picture below where the microphone is represented by system $\mathfrak H$:

\begin{center}
\begin{equation}\label{system}
\setlength{\unitlength}{.2cm}
\begin{picture}(43,8)(0,5)
\put(16.5,10){\makebox(0,0){$x(t)$}}
\put(31,10){\makebox(0,0){$y(t)$}}
\put(12.5,8){\vector(1,0){7.5}}
\put(28,8){\vector(1,0){7.5}}
\put(20,6){\framebox(8,4){$\mathfrak H$}}
\end{picture}
\end{equation}
\end{center} 
In the second step, the analog signal $\{y(t); t\in T\}$  becomes, in turn, the input of another system called CODEC which converts $y$ into a digital signal, namely a sequence of discrete  values from a finite alphabet like, for instance, $\{0,1\}$. Several other elaborations, carried out by suitable systems, are needed before our friend can hear the message. Coding a vocal input into bits and back-recode it into an audible (and hopefully understandable!)  output is one of the many applications of signal processing that we use on a daily basis.

A class of very useful systems (\ref{system}) are those given by a {\em convolution} of the form
\begin{equation}\label{conv}
y(t)=\int_{-\infty}^\infty h(t-s)x(s)ds.
\end{equation}
Here the signal $h(t)$, called impulse response, fully characterizes the system. These systems satisfy two important properties: They are {\em linear} and {\em time-invariant}. Linearity amounts to the superposition property that if $x_1$ and $x_2$ are transformed into $y_1$ and $y_2$, respectively, then $c_1 x_1(t)+c_2x_2(t)$ is transformed into $c_1 y_1(t)+c_2 y_2(t)$. Here $c_1$ and $c_2$ are any two numbers and it is implicitly assumed that the filter acts on a family of signals which is closed under linear combination (vector space) so that $c_1 x_1(t)+c_2x_2(t)$ is still in the family. Time invariance amounts to the fact that any time shift in the input causes exactly the same time shift in the output. Formally, $\{x(t-T); t_1\le t\le t_2\}$ gets transformed into $\{y(t-T); t_1\le t\le t_2\}$, for any $x$ and for any $T$.

Define the complex exponential
\begin{equation}\label{compexp}e^{ik2\pi ft}=\cos(k2\pi ft)+i\sin(k2\pi ft)
\end{equation}
It follows from the fundamental identity
\[\cos^2(\alpha)+\sin^2(\alpha)=1
\]
that the values of a complex exponential are on the unit circle, namely
\[|e^{ik2\pi ft}|=\sqrt{\cos^2(k2\pi ft)+\sin^2(k2\pi ft)}=1.
\]
The system (\ref{conv}) is fully characterized also by its {\em frequency response}
\begin{equation}\label{FR}
H(i2\pi f)=\int_{-\infty}^\infty h(t)e^{-ik2\pi f t}dt.
\end{equation}
This is called the {\em Fourier transform} of the impulse response and may be defined in a stronger or weaker sense depending on the properties of $h$.

Consider a convolutional filter (\ref{conv}) with frequency response $H$ (\ref{FR}). Suppose that the filter is {\em BIBO-stable} (cf. Subsection III\ref{frequency selective}). Then, if the input is the complex exponential $x(t)=e^{ik2\pi f t}$ (\ref{compexp}),
the corresponding output is simply 
\[y(t)=H(ik2\pi f)e^{ik2\pi f t}.
\] 
The complex exponential signal is namely mapped into {\em itself} times the {\em constant} $H(ik2\pi f)$. Mathematically, one says that complex exponentials are {\em eigenfunctions} of such filters. This property has an important consequence. Suppose now that the input is
\[x(t)=\sum_{k=-\infty}^\infty a_k e^{ik2\pi ft}.
\]
(This is just a Fourier series representation alternative to (\ref{FSeries})). Assume that the signal has finite energy on the period $[0,1/f]$ which is equivalent to
\[\sum_{k=-\infty}^\infty |a_k|^2<\infty.
\]
(This happens if the {\em Fourier coefficients} $a_k$ decay sufficiently rapidly as $|k|\rightarrow\infty$). Then, thanks to linearity and stability of the filter plus the eigenfunction property, the signal $x$ gets mapped to 
\[y(t)=\sum_{k=-\infty}^\infty b_k e^{ik2\pi ft}
\]
where
\begin{equation}\label{selective}
b_k=H(ik2\pi f)\cdot a_k.
\end{equation}
This implies the relation between squared absolute values
\begin{equation}\label{ABSselective}|b_k|^2=|H(ik2\pi f)|^2\cdot |a_k|^2
\end{equation}
describing how the average power (energy rate) of the $k$-th harmonic has changed. This result is of paramount importance. For instance, if we design the filter so that its frequency response has small absolute value in a certain frequency band $[f_1,f_2]$, then all harmonics of the output with $kf\in[f_1,f_2]$ will be attenuated. In other words, (\ref{selective}) opens the way to {\em selective manipulation of the harmonics} !

\subsubsection{Low-pass vs high-pass filters}\label{low vs. high filters}
Consider a damped hamonic oscillator subject to an external force $F(t)=x(t)$. This may represent a mass-spring system in a fluid and subject to gravity. Let $y(t)$ represent the displacement of the mass $m$ barycenter from its equilibrium position. Then, by Newton's law, $y$ obeys the differential equation
\begin{equation}\label{Newton}m\frac{d^2y}{dt^2}+\beta\frac{dy}{dt}+Ky=x, \quad m>0, \quad \beta>0, \quad K>0.
\end{equation}
The  frequency response $H$ (\ref{FR}) of this input-output system transforming the external force into the position of the barycenter is readily seen to be
\[H(i\omega)=\frac{1}{m(i\omega)^2 + \beta i \omega+K}, \quad \omega=2\pi f.
\]
As the circular frequency $\omega$ or, equivalently, the frequency $f$ increases to $+\infty$, $|H(i\omega|$ decreases to zero. Thus this system can be considered as an approximation of an ideal low-pass filter which considerably attenuates the high-frequency harmonics in view of (\ref{ABSselective}).  Moreover, and differently form the ideal filter, this filter is {\em stable}.  Stability here is due to the presence of the friction force $-\beta\frac{dy}{dt}$ causing the damping. Indeed, an undamped harmonic oscillator is an unstable system which may transform bounded signals into unbounded ones as in the resonance phenomenon well-known, for instance, in acustics. 

\begin{figure}
\begin{center}
	\includegraphics[scale=1.5]{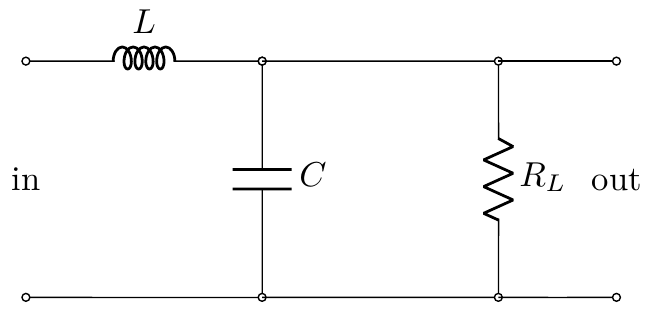}
	\caption{A RLC circuit acting as a low-pass filter.}
	\label{RLC}
	\end{center}
\end{figure}

Mechanical harmonic oscillators and their electric network counterpart (RLC circuits) play a central role in physics and in many  artificial devices. 
In Figure \ref{RLC},
a second-order {\em passive} filter is depicted. This RLC (Resistor-Inductor-Capacitor) circuit is the analog counterpart of the mechanical system (\ref{Newton}). It is named passive because it is made only from passive components,  namely it does not require an external power source beyond the signal. Active filters, on the contrary, feature one or more {\em active} components like the {\em operational amplifier} of Figure \ref{active} which does require an external power source.

Consider now the ideal differentiator of Subsection 2.1 $y(t)=\frac{dx}{dt}(t)$. 
This filter  is unstable but it is often used in control and telecommunications systems because of the positive characteristics described in Subsection III\ref{frequency selective}. Its frequency response (\ref{FR}) is $H(i\omega)=i\omega, \omega=2\pi f$, so  that the magnitude is $|H(i\omega)|=|\omega|$. Thus, this filter increases the energy of the high-frequency harmonics {\em proportionally to their frequency}.   Input-output stability, which is lost when a differentiator is present in the open-loop setting (\ref{system}), is typically recovered in information engineering using {\em feedback}
\begin{center}
$$
\setlength{\unitlength}{.2cm}
\begin{picture}(30,14)(0,0)
\put(3,11.5){\makebox(0,0){$r$}}
\put(6.5,11.5){\makebox(0,0){$+$}}
\put(10.5,11.5){\makebox(0,0){$x$}}
\put(29.5,11.5){\makebox(0,0){$y$}}
\put(9.5,8.5){\makebox(0,0){$-$}}
\put(0,10){\vector(1,0){7}}
\put(9,10){\vector(1,0){3}}
\put(26,10){\vector(1,0){8}}
\put(8,2){\vector(0,1){7}}
\put(8,10){\circle{2}}
\put(28,10){\circle*{0.5}}
\put(28,2){\line(0,1){8}}
\put(8,2){\line(1,0){9.5}}
\put(28,2){\vector(-1,0){11}}
\put(12,8){\framebox(14,5){$H(i\omega)=i\omega$}}
\end{picture}
$$
\end{center}

In electronics, a simple first-order passive RC filter may be used to approximate the ideal differentiator
\begin{equation}\label{RC}RC\frac{d V_{{\rm out}}}{dt}(t) +V_{{\rm out}}(t)= RC\frac{d V_{{\rm in}}}{dt}(t).
\end{equation}
The approximation, however, is  only good for low frequencies. A more precise differentiator may be realized using an {\em active} element such as an {\em operational amplifier}, cf. Figure \ref{active}.

\subsection{Aliasing}\label{aliasing}

An anti-aliasing filter is a filter that removes the high-frequency of a signal to be {\em sampled}. This, in order to avoid producing artifacts (alias = non existing low-frequency harmonics), due to undersampling, in the reconstruction of the signal, see \cite[Section 7.3]{oppenheim1997signals}. A well-known example of aliasing is that in movies we sometimes see helicopter blades or stagecoach wheels rotating slower than they actually are or even rotating in the opposite direction. The sampling is here due to the frames, e.g. 24 per second, of the movie.
As a trivial example of aliasing, consider the $k=5$ normal mode  $\sin(5\cdot \pi t)=\sin(5\pi t)$ of the vibrating string of Subsection II\ref{harmonics} (Figure \ref{vibratingstring}). Suppose samples of the vertical displacement of the string are taken at positions $0, 1/5, 2/5, 3/5, 4/5$ where such displacement  is always zero ($\sin(0)=\sin(2\pi)=\sin(4\pi)=\sin(6\pi)=\sin(8\pi)=0$). Then, the reconstructed signal will be the first harmonic identically equal to zero (straight line) corresponding to $k=0$. Thus, due to undersampling, the $k=5$ harmonic has been replaced by the $k=0$ harmonic and the circular frequency  $\omega_5=5\cdot \pi$ has taken the identity  of the lower frequency $\omega_0=0\cdot \pi=0$ (aliasing). Suppose now samples are taken at positions $0,1/3, 2/3$. Then the vertical displacements are $0,-\frac{\sqrt{3}}{2},-\frac{\sqrt{3}}{2}$. The reconstructed signal can be shown to be  $\sin(-\pi t)$ corresponding to $k=-1$ (depicted together with the $k=1$ harmonic at the top of Figure \ref{vibratingstring}) and, again, we have the aliasing phenomenon. By Shannon's interpolation theorem \cite[Chapter 7]{oppenheim1997signals}, if we sample $\sin(10\pi t)$ rather than at multiples of $1/5$ or $1/3$ at multiples of any $T<\frac{1}{10\pi}$, we can reconstruct the signal exactly.

\subsection{Historical bits}\label{bits}

\noindent
{\em Harmonics}

The concept of {\em harmony} was central in the Pythagorean school of philosophy in the sixth century B.C.E. in Crotone, Southern Italy. While harmony, according to this school, pervades the universe, it is based on the divinity of the natural numbers $1, 2, 3,\ldots$. Pythagoreans, in their experimental studies on a chord musical  instrument, discovered    {\em overtones} which had a frequency that was an integer multiple of a fundamental frequency, i.e. harmonics.

\vspace{.5cm}
\noindent
{\em Fourier}

On December 21st, 1807, Joseph B.\ Fourier submits a manuscript to the Institute of France in Paris entitled  “Sur la propagation de la chaleur”. He had been working on this memoir during his stay in Grenoble as Prefect of the Department of Is\'ere, a post to which he had been appointed by Napoleon himself. The surprising, but not fully justified, results provoke an animated discussion among the examiners. The committee consists of Lagrange, Laplace, Lacroix and  Monge. Lagrange and Laplace, who criticize Fourier's expansion  of functions as trigonometrical series, and Monge had all been teachers of Fourier at the \'Ecole Normale.  Moreover, Fourier and Monge had also been together in the expedition to Egypt in 1798 and members of the mathematics division of the Cairo Institute together with Malus and Napoleon himself. In spite of all of this,  the manuscript is eventually rejected, showing how controversial the idea was initially. Fourier's  “Th\'eorie analytique de la chaleur”  will appear only in 1822.

\vspace{.5cm}
\noindent
{\em Feedback}

Feedback was employed in ancient water clocks by the Greeks and the Arabs and, extensively, in mechanical clocks and other devices starting from the Middle Ages. In 1788, James Watt designed his centrifugal governor to automatically regulate the angular velocity of the steam engine. The first mathematical paper on feedback was published by James Clerk Maxwell (the father of electromagnetism) only in 1868 \cite{maxwell1868governors}. The early basis of the modern theory was laid down around 1930 by Blake, Bode and Nyquist who were interested chiefly in constructing amplifiers for telephone communication which would attenuate noise. Important developments took place during the space race \cite{kalman1960new} and in the seventies, see e.g. \cite{doyle2013feedback}.


\section*{Acknowledgments}
 Authors are listed in alphabetical order as they  contributed equally to the manuscript.
 We are thankful to the Heuss scholarship, to the New School for Social Research, and to New York for providing the common ground, literally and metaphorically, for the idea of the present manuscript. We also thank  Gabriele Vianelli for his comments on harmonics in music. We are grateful to Daniele Alpago for preparing the figures.

\bibliographystyle{apacite}
\bibliography{refs.bib}

\end{document}